\documentclass[12pt]{article}
\usepackage{amsfonts,amsmath,amsthm,graphicx}
\newcommand{\p}{{\mathbb P}}
\newcommand{\R}{\mathbb{R}}
\newcommand{\E}{\mathbb{E}}
\newcommand{\Z}{\mathbb{Z}}
\newcommand{\var}{\text{Var}}
\newcommand{\vol}{\text{Vol}}
\newcommand{\pd}{U^{\text{diff}}}

\newtheorem{thm}{Theorem}
\newtheorem{prop}[thm]{Proposition}
\newtheorem{lem}[thm]{Lemma}
\newtheorem{cor}[thm]{Corollary}
\newtheorem{defi}[thm]{Definition}

\def\R{\mathbb{R}}

\def\P{\mathbb{P}}
\def\E{\mathbb{E}}
\def\N{\mathbb{N}}
\def\Z{\mathbb{Z}}

\def\eps{\varepsilon}

\title{Gravitational allocation to Poisson points}
\author{Sourav Chatterjee\thanks{U.C. Berkeley. Research
supported in part by NSF grants \#DMS-0244479 and \#DMS-0104073.}
 \and Ron Peled$^*$
  \and Yuval Peres\thanks{Microsoft Research and U.C. Berkeley. Research
supported in part by NSF grants \#DMS-0244479 and \#DMS-0104073.}
 \and
Dan Romik\thanks{Bell Laboratories. Research
supported in part by NSF grants \#DMS-0244479 and \#DMS-0104073.}
}
\begin{document}
\maketitle

\vspace{-35.0pt}
\begin{abstract}
For $d\ge 3$, we construct a non-randomized, fair and translation-equivariant allocation of Lebesgue measure to the points of a standard Poisson point process in $\R^d$, defined by allocating to each of the Poisson points its basin of attraction with respect to the flow induced by a gravitational force field exerted by the points of the Poisson process. We prove that this allocation rule is economical in the sense that the \emph{allocation diameter}, defined as the diameter $X$ of the basin of attraction containing the origin, is a random variable with a rapidly decaying tail. Specifically, we have the tail bound
$$ \p(X > R) \le C \exp\Big[-c R (\log R)^{\alpha_d} \Big]$$
for all $R>2$, where: $\alpha_d = \frac{d-2}{d}$ for $d\ge 4$; $\alpha_3$ can be taken as any number $<-4/3$; and $C,c$ are positive constants that depend on $d$ and $\alpha_d$.
This is the first construction of an allocation rule of Lebesgue measure to a Poisson point process with subpolynomial decay of the tail $\p(X>R)$.
\end{abstract}

\bigskip \noindent AMS 2000 subject classification: 60D05.

\bigskip \noindent Key words: Gravitation, fair allocation, Poisson process, translation-equivariant mapping.


\newpage

\tableofcontents

\newpage

\section{Introduction \label{sectionintro}}

Let $d \in \mathbb{N}$. Let $\Xi$ be a discrete subset of $\R^d$. We call the elements of $\Xi$
{\bf centers}.
An {\bf allocation} (of Lebesgue measure to $\Xi$) is a measurable function $\psi:\R^d\to \Xi\cup\{\infty\} $ that satisfies
\begin{eqnarray*}
\vol(\psi^{-1}(\infty)) &=& 0, \\
\vol(\psi^{-1}(z)) &=& 1, \qquad z\in \Xi,
\end{eqnarray*}
where $\vol(\,\cdot\,)$ is Lebesgue measure in $\R^d$. For $z\in \Xi$, we call $\psi^{-1}(z)$ the {\bf cell allocated to $z$}. In other words, an allocation is a way of partitioning $\R^d$ into cells of Lebesgue measure 1 that together cover $\R^d$ up to a set of measure $0$, and assigning them to the points of $\Xi$.

Let $Z$ be a translation-invariant simple point process in $\R^d$ with unit intensity defined on some probability space $(\Omega,{\cal F},\p)$. That is, $Z$ is a random discrete subset of $\R^d$ such that for
any open set $A\subset \R^d$, the random variable $|A\cap Z|$ (where $|E|$ denotes the cardinality of a set $E$) has mean $\vol(A)$, and for
any $x\in\R^d$ and open sets $A_1,A_2,\ldots,A_k\subset \R^d$, the random vector $(|(A_1+x)\cap Z|,|(A_2+x)\cap Z|,\ldots,|(A_k+x)\cap Z|)$ has distribution that does not depend on $x$.
An {\bf allocation rule} (of Lebesgue measure to $Z$) is a  mapping $Z\to \psi_Z$ that is defined $\p$-a.s., measurable (with respect to the relevant $\sigma$-algebras), and such that: (i) almost surely $\psi_Z$ is an allocation of Lebesgue measure to $Z$, and (ii) the mapping $Z\to \psi_Z$ is translation-equivariant, in the sense that $\p$-a.s., for any $x,y\in\R^d$ we have
$$ \psi_{Z+x}(y+x) = \psi_Z(y) + x. $$

\begin{figure}[h!]
 \hspace{-40.0pt} \resizebox{!}{14.0cm}{\includegraphics{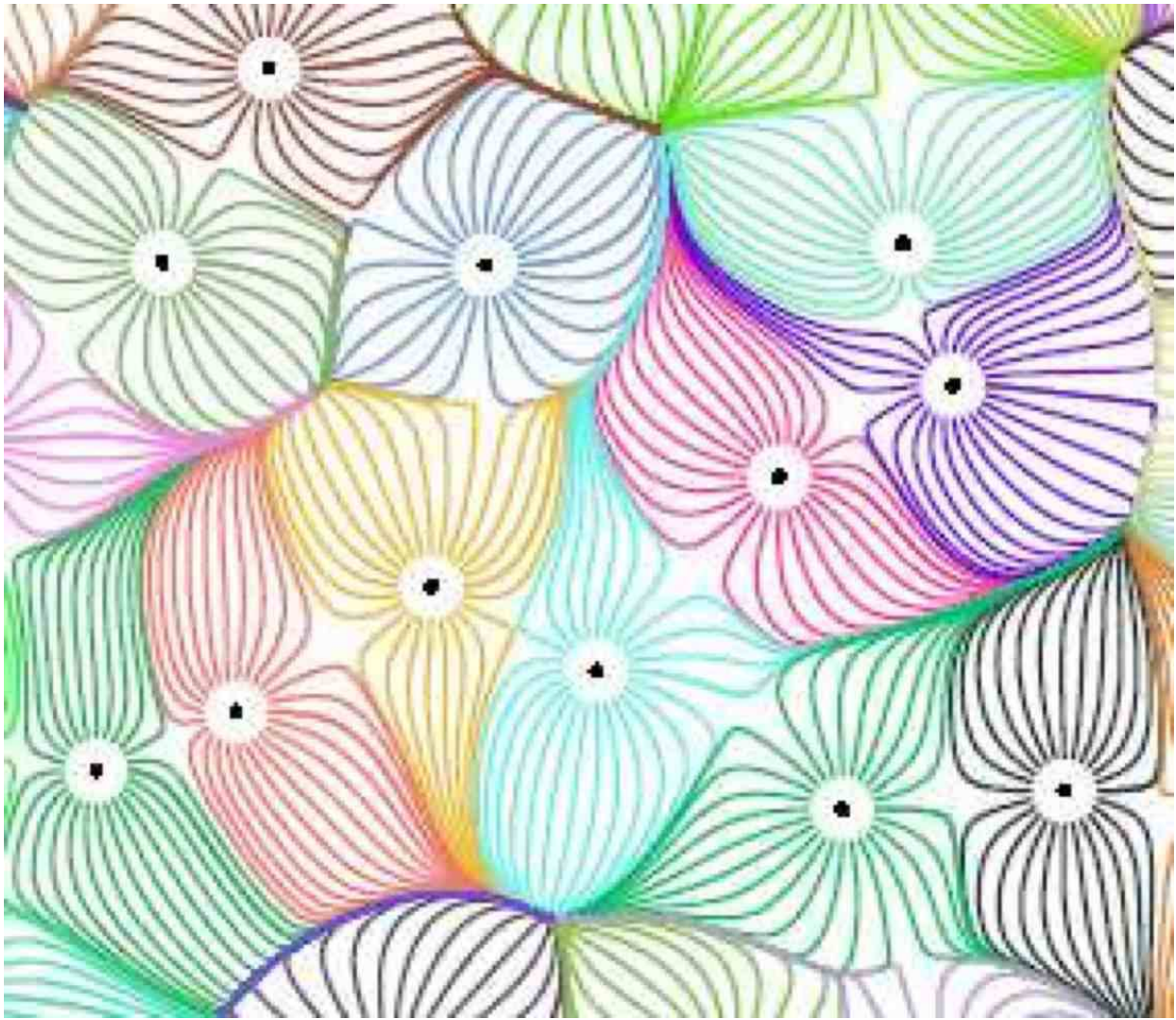}}
\caption{The Nazarov-Sodin-Tsirelson-Volberg gradient flow allocation, equivalent to 2-dimensional gravitational allocation. It can be defined for a finite point set or for the process of zeros of the Gaussian Entire Function. For the Poisson point process we construct the analogous allocation  in dimensions 3 and higher.
(Picture due to Manjunath Krishnapur). \label{manjupic}}
\end{figure}

Figure \ref{manjupic} shows a particularly important example of an allocation rule that gave much of the inspiration for the current paper -- see below.

An allocation rule $Z\to \psi_Z$ may satisfy several additional desirable properties: each cell
$\psi_Z^{-1}(z)$ may be open,
may contain its ``owner'' $z$; each cell $\psi_Z^{-1}(z)$ may be connected; each cell may be bounded. In the event that a.s.\ all the cells are bounded, one may consider the {\bf allocation diameter}, which is the random variable
$$ X = \text{diam}(\psi_Z^{-1}(\psi_Z(0))), $$
where $\text{diam}(\,\cdot\,)$ denotes the diameter of a set. The rate of decay of the tail $\p(X>R)$ of the distribution of $X$ as $R\to\infty$ can be used as a quantitative measure for how economical the allocation rule is; roughly, a fast rate of decay means that it is rarer for points to be allocated to a far-away location. Note that by translation-equivariance one may take the diameter of the cell $\psi_Z^{-1}(\psi_Z(x))$ containing any given point $x\in\R^d$ and get a random variable equal in distribution to $X$.

Holroyd and Peres \cite{holroydperes} showed that if $d=1,2$ and $Z$ is a standard Poisson point process of unit intensity in $\R^d$, then for every allocation rule the allocation diameter $X$ satisfies $\E X^{d/2} = \infty$. In particular, in this case the decay of $\P(X>R)$ to 0 cannot be faster than polynomial in $R$.
Hoffman, Holroyd and Peres \cite{stablemarriage1} constructed an allocation rule for every translation-invariant point proces in $\R^d$ with unit intensity, the {\bf stable marriage allocation}, in every dimension $d\ge 1$. In the stable marriage allocation, almost surely
all the cells are open, bounded and contain their owners, but not all are connected, and when $Z$ is a Poisson point process the allocation diameter $X$ satisfies $\E X^d=\infty$. The stable marriage allocation rule is the unique one which has the ``stability'' property that almost surely, for any $z,z'\in Z$ and $x,y\in \R^d$, if $\psi_Z(x) = z$ and $\psi_Z(y)=z'$
then $|x-z'|\ge |x-z|\wedge |y-z'|$.

Nazarov, Sodin and Volberg \cite{nazarovetal2} recently constructed an allocation rule based on an idea suggested by Tsirelson in \cite{sodintsirelson}. Their allocation rule is defined for the two-dimensional point process ${\cal X}$ of zeros of the Gaussian Entire Function (GEF), which is the random analytic function
$$ f(z) = \sum_{n=0}^\infty \xi_n \frac{z^n}{\sqrt{n!}}, \qquad z\in \mathbb{C}, $$
where $(\xi_n)_{n=0}^\infty$ are i.i.d.\ standard complex gaussian random variables. In their construction, the cell of each $z\in {\cal X}$ is defined as the {\bf basin of attraction} of $z$ with respect to the flow induced by the random planar vector field $z\to \big(\nabla \log|f|\big)(z)- z$. The cells are
connected by definition, and in \cite{nazarovetal2} it was proved
that they are a.s.\ bounded, each have area $\pi$ (which is the
reciprocal of the mean density
of points in the process of zeros of the GEF),
and that there exist absolute constants $C,c>0$ such that
the allocation diameter $X$ satisfies
$$ c e^{-C R(\log R)^{3/2}} \le \p(X\ge R) \le C e^{-c R (\log R)^{3/2}}, \qquad R>1. $$

\begin{figure}[h!]
\centering \resizebox{!}{9.0cm}{\includegraphics{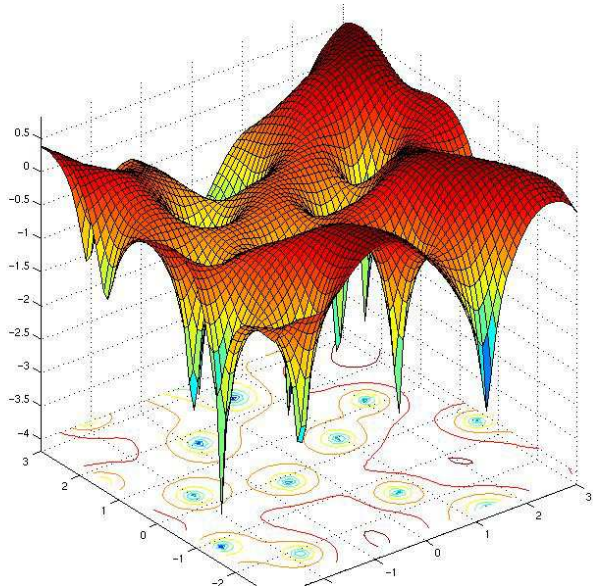}}
\caption{The potential function associated with planar gravitational
allocation
(picture due to Manjunath Krishnapur).\label{potentialgraph}}
\end{figure}

\begin{figure}[h!]
\centering \resizebox{!}{9.0cm}{\includegraphics{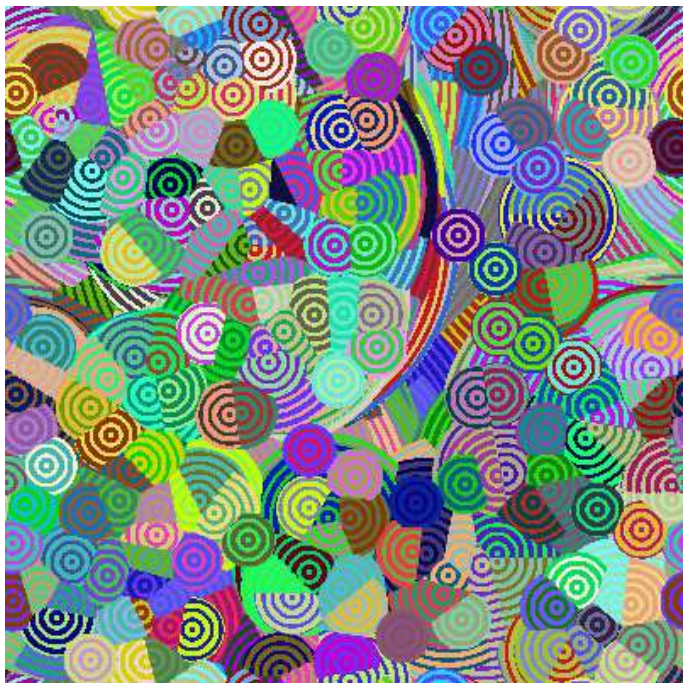}}
\caption{The two-dimensional stable marriage allocation for a Poisson process (picture due to Alexander E. Holroyd).\label{picturesalloc}}
\end{figure}

Figure \ref{manjupic} shows a simulation of the
Nazarov-Sodin-Tsirelson-Volberg {\bf gradient flow allocation}.
Figure \ref{potentialgraph} shows the graph of the potential function
$\log |f|$ associated with the allocation (where $f$ is an
approximation to the GEF).
Figure \ref{picturesalloc} shows a simulation of the stable marriage
allocation in 2 dimensions.

In this paper, we construct a new allocation rule of Lebesgue measure to the points of the standard Poisson point process in $\R^d$, for any $d\ge 3$. Our construction was inspired by the gradient flow allocation, and we call it {\bf gravitational allocation}. To define it, denote by ${\cal Z}$ the standard Poisson process in $\R^d$. Consider
the random vector field $F:\mathbb{R}^d \to \mathbb{R}^d$ defined by
\begin{equation} \label{eq:forcefield}
F(x) = \sum_{z\in {\cal Z},\ |z-x|\uparrow} \frac{z-x}{|z-x|^d},
\end{equation}
where the summands are arranged in order of increasing distance from $x$. The term
$(z-x)/|z-x|^d$ represents a gravitational force felt by a unit mass at a point $x$ due to the influence of a unit mass placed at point $z$. When $d=3$, this is the ordinary Newtonian gravitational force.
An elementary observation that can be traced back to Chandrasekhar \cite{chandrasekhar} (see also \cite{heathshepp}), based on the Kolmogorov three-series theorem, is that for any fixed $x\in\R^d$, the infinite series for $F(x)$ converges almost surely (this is true for $d\ge 3$). The random vector $F(x)$ has a symmetric stable distribution of index $\frac{d}{d-1}$. This can be seen using a simple scaling argument (see Remark (iv) below), or by an exact computation, see \cite{heathshepp}.

We prove the following result concerning the process of gravitational forces acting simultaneously on all points of $\R^d$:

\begin{prop}[Simultaneous\ \ convergence\ \ and\ \ differentiability] \label{prop1} Assume $d\ge 3$.
Almost surely, the series in \eqref{eq:forcefield} converges simultaneously for all $x$ for which it is defined (namely, all $x\in \R^d \setminus {\cal Z}$) and defines a translation-invariant (in distribution) vector-valued random function. The random function $F$ is almost surely continuously differentiable where it is defined.
\end{prop}

Note that since the sum in \eqref{eq:forcefield} does not converge absolutely, the choice of the order of summation is essential for Proposition \ref{prop1} and the results below to hold.

Consider now the integral curves $Y(t)$ of the vector field $F$, that is, solutions of the equation
$$ \dot{Y}(t) = F(Y(t)). $$
We call these curves the {\bf gravitational flow curves} (in a simplified inertia-less Newtonian gravitational world). For $x\in \R^d\setminus {\cal Z}$, denote by $Y_x$ the integral curve with initial condition $Y_x(0)=x$. By Proposition \ref{prop1} and standard ODE existence and uniqueness theorems, $Y_x$ is defined up to some maximal positive time $\tau_x$ (where possibly $\tau_x=\infty$). For each center $z\in {\cal Z}$, say that the curve $Y_x$ {\bf ends at $z$} if $\lim_{t\uparrow \tau_x} Y_x(t) = z$, and define the {\bf basin of attraction} of $z$ by
$$ B(z) = \{ x\in \R^d \setminus{\cal Z}\ |\ Y_x(t)\text{ ends at }z\}\cup \{ z \}. $$
Define $$\psi_{\cal Z}(x) = \left\{ \begin{array}{ll} z\qquad & x\in B(z)\text{ for }z\in {\cal Z},\\ \infty &
x\notin \bigcup_{z\in {\cal Z}} B(z). \end{array}\right. $$
Our main result is the following theorem.

\begin{thm}[Fairness and efficiency of the allocation] \label{thm2}
The mapping ${\cal Z}\to\psi_{\cal Z}$ is an allocation rule of
Lebesgue measure to the Poisson point process ${\cal Z}$. Almost
surely all the cells $\psi_{\cal Z}^{-1}(z)$ are bounded.
The allocation diameter $X = \text{diam}(\psi_{\cal Z}^{-1}(\psi_{\cal Z}(0)))$ satisfies the following tail bounds:
In dimensions 4 and higher, we have
\begin{equation} \label{eq:tailestimate1}
\p(X > R) \le C_1 \exp\Big[-c_2 R (\log R)^{\frac{d-2}{d}} \Big]
\end{equation}
for some constants $C_1=C_1(d),c_2=c_2(d)>0$ and all $R>2$.
In dimension 3, for any $\alpha>0$ there exist constants $C_1=C_1(\alpha),c_2=c_2(\alpha)>0$ (depending on $\alpha$) such that for all $R>2$ we have
\begin{equation}\label{eq:tailestimate2}
\p(X > R) \le C_1 \exp\bigg[-c_2 \frac{R}{(\log R)^{\frac{4}{3}+\alpha}} \bigg].
\end{equation}
\end{thm}

\begin{figure}[h!]
\begin{center}
\centering \resizebox{!}{10.0cm}{\includegraphics{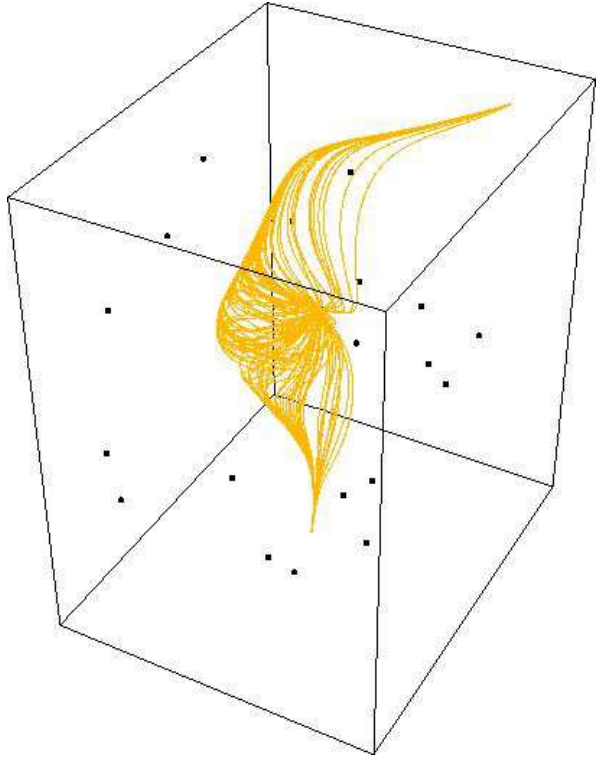}}
\caption{Simulation of a cell in 3-dimensional gravitational allocation \label{gravalloc}}
\end{center}
\end{figure}

Note that the cells in gravitational allocation are open, connected and contain their owners. They are also contractible, see Remark (v) below. Figure \ref{gravalloc} shows a simulation of a cell in 3-dimensional gravitational allocation.

For $L>0$ and $x\in\R^d$ denote by $Q(x,L)$ the box $x+[-L,L]^d$.
A main ingredient in the proof of Theorem 2 is the following result.

\begin{thm}[Bounds for the probability of an $R$-crossing] \label{thm3}
Let $E_R$ denote the event that there exists an integral curve $Y(t)$ connecting $\partial Q(0,R)$ and $\partial Q(0,2R)$ (in either order). Then, if $d\ge 4$ then we have
$$
\p(E_R) \le C_1 \exp\Big[-c_2 R (\log R)^{\frac{d-2}{d}} \Big]$$
for some constants $C_1=C_1(d),c_2=c_2(d)>0$ and all $R>2$.
In dimension 3, for any $\alpha>0$ there exist constants $C_1=C_1(\alpha),c_2=c_2(\alpha)>0$ (depending on $\alpha$) such that for all $R>2$ we have
$$ \p(E_R) \le C_1 \exp\bigg[-c_2 \frac{R}{(\log R)^{\frac{4}{3}+\alpha}} \bigg].$$
\end{thm}

In a forthcoming paper \cite{gravpaper2}, we will prove
lower bounds for the tail of the distribution of the allocation diameter $X$,
and additional bounds on the distance $|\psi_{\cal Z}(0)|$ of the origin from
its star.

\paragraph{Further remarks.}
{\bf (i)} Allocation rules have an equivalent description as non-randomized {\bf extra head rules}. If $Z$ is a translation-invariant simple point process of unit intensity in $\R^d$, an extra head rule for $Z$ is a random variable $T$ coupled with $Z$ so that a.s.\ $T\in Z$ and the random set $Z-T$ has the same distribution as $Z$ conditioned to have a point at $0$. The extra head rule is said to be non-randomized if $T$ is measurable with respect to $Z$. In \cite{holroydperes} it was shown that if $\psi_Z$ is an allocation rule then $T=T_Z=\psi_Z(0)$ is a non-randomized extra head rule, and conversely,  given a non-randomized extra head rule $T_Z$, the mapping $\psi_Z(x) = T_{Z-x}$ is an allocation rule.
\\
{\bf (ii)} For any $u\in \R^d$ one may replace the vector field $F(x)$ by $F(x)+u$ and obtain a modified allocation rule. Thus, there is more than one possible construction of an allocation rule involving the gravitational field, and one might speculate that a suitable modification of the construction might lead to better tail bounds for the allocation diameter.
\\
{\bf (iii)} For some results on the related topic of translation-invariant perfect matchings for point processes, see \cite{holroydperesmatchings}. For related results on matchings between random point configurations in a finite setting, see the papers \cite{ajtaietal, leightonshor, talagrand}.
\\
{\bf (iv)} Here is a simple argument proving that for fixed $x$ the force vector $F(x)$ has a stable distribution with scaling exponent $d/(d-1)$ (we believe this argument is known but could not locate a reference; the proof of this fact in \cite{heathshepp} uses explicit computations and is more complicated). If $F_1, F_2, \ldots, F_n$ are i.i.d. copies of $F(x)$, then their sum is the force exerted on $x$ by the union of $n$ independent copies of the Poisson process, which is a Poisson process with intensity $n$ (or equivalently a Poisson process of unit intensity scaled by $n^{-1/d}$). Thus, because the individual force terms scale as the $(d-1)$th power of the distance, by rescaling it follows that $F_1+F_2+\ldots+F_n$ has the same law as $n^{(d-1)/d} F(x)$, which proves our claim.
\\
{\bf (v)} Another interesting property of gravitational allocation is that the cells are contractible. This is immediate from their definition as the basins of attraction with respect to the flow of the vector field $F$. Formally, denote by $(\Phi_t)_{t\ge 0}$ the flow semigroup of the vector field $F$, and for each $x\in \R^d$ denote by $\tau_x$ the time for $x$ to flow to its star $\psi_{\cal Z}(x)$ (that is, the maximal time for which the curve $Y_x$ is defined, or $0$ for the star). Then if $z\in {\cal Z}$ and $B(z)$ is its basin of attraction, the mapping $\varphi:B(z)\times [0,1]\to B(z)$ defined by $\varphi(x,t)= \Phi_{\tau_x t}(x)$ is a homotopy between the identity map $\text{id}_{B(z)}$ and the constant mapping $B(z)\to z$. (Note that almost surely, for all $x\in B(z)\setminus \{z\}$ we have $\tau_x<\infty$, since by definition we have that $\lim_{t\uparrow \tau_x} Y_x(t)=z$, and $F(u)=(z-u)/|z-u|^d+O(1)$ when $u\to z$, so that once the flow curve $Y_x(t)$ approaches $z$, it must reach $z$ in a finite time.)
\bigskip

\paragraph{A reading guide.}
Here is a guide to reading the rest of the paper.
Section \ref{sectionnotation} introduces some notation and recalls some standard
estimates for Poisson random variables.
In section \ref{sectionoutlines} we give outlines of the
proofs of the main claims, which we hope will give the reader a higher-level
picture of the ideas in the paper and will simplify reading the more technical
later sections. In section \ref{sectionproofthm2} we show how the main result,
Theorem \ref{thm2}, can be deduced fairly easily from Theorem 3 (bounds for
the probability of an $R$-crossing). These sections are easy to read and
we recommend starting with them.

The remaining sections constitute the main technical parts of the paper.
Sections \ref{sectionproofprop1}, \ref{sectionforcethm}, \ref{sectionpotential}
and \ref{sectionlarge} are the ``preparation'' part:
In Section \ref{sectionproofprop1} we prove Proposition \ref{prop1}
(simultaneous convergence and differentiability of the force)
and Proposition \ref{rearrangement} (a useful alternative formula for
the force, see Section \ref{sectionoutlines} below). Section
\ref{sectionpotential} contains a similar but slightly more difficult analysis
for the {\bf gravitational potential} function, an auxiliary function that
is defined only in dimensions 5 and higher. In Section \ref{sectionforcethm}
we prove an important auxiliary theorem bounding the joint density
of a vector of forces. In Section \ref{sectionlarge} we prove large deviations
results that will be used repeatedly as the main ``engine'' in the proof of
Theorem \ref{thm3}.

Finally, Sections \ref{sectionproofthm3} and \ref{sectiondim34} contain
the proof of Theorem 3. The proof is split into two parts. The first and
slightly simpler case is the
proof in dimensions 5 and higher. The last section, Section \ref{sectiondim34},
treats the more delicate case of dimensions 3 and 4.
We recommend to the reader who is mainly interested in our main result
to only skim through the results in Sections
\ref{sectionproofprop1}, \ref{sectionforcethm}, \ref{sectionpotential}
and \ref{sectionlarge} and to proceed to Section
\ref{sectionproofthm3}. However,
we believe the results in these auxiliary sections to be of significant
independent interest.

\paragraph{Acknowledgements}
We are grateful to Manjunath Krishnapur, Misha Sodin and Fedor Nazarov for helpful comments and discussions, and to an anonymous referee for suggesting innumerable valuable comments and corrections, including a shorter and more elegant proof of the uniqueness of solution of a system of equations that forms part of the proof of Theorem \ref{joint_force_density_thm}.

\section{Preliminaries \label{sectionnotation}}

Here is some notation that we will use throughout the paper: $d$ denotes the dimension, and will always be an integer $\ge 3$ (in some theorems it will be assumed explicitly that $d\ge 5$ or that $d=4$ or that $d=3$).
We denote by $|x|$ the Euclidean norm of a vector. We denote Lebesgue measure in $\R^d$ by $\vol(\,\cdot\,)$ . If $V=(V_1,...,V_k)$ is a random vector, we denote by $\var(V)$ the sum of the variances of its coordinates. Let $\kappa_d = \pi^{d/2}/\Gamma(d/2+1)$ be the volume of the unit ball in $\R^d$. Denote by ${\cal Z}$ the Poisson process of unit intensity on $\R^d$, and by $\p$ the probability measure on the probability space on which it is defined. For concreteness we denote ${\cal Z}=(z_i)_{i=1}^\infty$ for the specific ordering of the points of ${\cal Z}$ by increasing distance from $0$. For obvious reasons we refer to the $z_i$ as {\bf stars}. The letters $C,c$ (possibly with subscripts) will be used to denote various positive constants that may depend on the dimension $d$, where $C$ will typically be a large positive constant and $c$ will typically be a small positive constant, and the same symbols (such as $c_1$, etc.) may be used in different places with different numerical values. Big-O notation will be used, and it is understood that all constants implicit therein may depend on $d$ (and occasionally on other parameters that are kept constant throughout the discussion). We denote by $B(x,L)$ the ball of radius $L$ around $x\in\R^d$, and by $Q(x,L)$ the box $x+[-L,L]^d$. This notation and other notations that are used frequently in the paper are summarized in Table 1, which may be used for reference.

\begin{center}
\begin{table}
\caption{Summary of the main notation used in the paper}
\vspace{15.0pt}
\begin{tabular}{l|l|p{270.0pt}}
Symbol & Sections & Meaning \\ \hline
$d$ & all & The dimension, an integer $\ge 3$. In Section 5, $d \ge 5$. \\
$B(x,t)$ & all & Ball of radius $t$ around $x\in \R^d$. \\
$Q(x,t)$ & all & The cube $x + [-t,t]^d$. \\
$\kappa_d$ & all & $\frac{\pi^{d/2}}{\Gamma(d/2+1)} = $volume of the unit ball in $\R^d$. \\
${\cal Z}=(z_i)_{i=1}^\infty$ & all & The ``stars'': a standard Poisson point process  in $\R^d$. \\
$F(x)$ & all & The random gravitational force field induced by ${\cal Z}$. \\
$E_R$ & 1, 3, 9 & The event of a gravitational flow curve crossing between $\partial Q(0,R)$
and $\partial Q(0,2R)$. \\
$g(x)$ & 5, 6, 8 & $g(x)=\frac{x}{|x|^d}$. \\
$D_k[\,\cdot\,]$ & 5, 6, 7, 8, 9 & The $k$-th derivative tensor of a function.\\
$R$ & all & The main parameter. \\
$B$ & 9, 10 & $R^{8/9}$ (In Section \ref{subsectiondim3}: $\frac{R}{(\log R)^\beta}$). \\
$\Delta$ & 9, 10 & A large constant. \\
$r$ & 9, 10& $\Delta\cdot (\log R)^{2/d}$ (In Section \ref{subsectiondim3}: $(\log R)^{1/3}\log \log R$). \\
$\rho$ & 9, 10& $R^{-1/10}$ (In Section \ref{subsectiondim3}: $\frac{1}{(\log R)^\gamma}$). \\
$s$ & 9, 10& $R^{-\frac{1}{10(d^2+1)}}$ (In Section \ref{subsectiondim3}: $\frac{1}{(\log R)^\delta}$). \\
$\varepsilon$ & 9, 10& $\frac{\rho}{s^d}\log R$ (In Section \ref{subsectiondim3}: $\frac{\rho}{s^3}$). \\
$\lambda$ & 9, 10 & $2\left(\frac{d}{\kappa_d}\right)^{1/d}(\log R)^{1/d}$
(In Section \ref{subsectiondim3}: $\sqrt{\log \log R}$). \\
$U(x)$ & 7, 8, 9 & The stationary centered gravitational potential. \\
$U(x\ |\ A)$ & 7, 8, 9, 10 & Centered contribution to the potential from stars in the set $A$. \\
$F(x\ |\ A)$ & 8, 9, 10 & Contribution to the force from stars in $A$. \\
$S$ & 9, 10 & The grid $S = r\Z \cap (Q(0,2R)\setminus Q(0,R))$. \\
$S_w$ & 9, 10 & For $w\in S$, the subgrid $s\Z \cap (Q(w,2r)\setminus Q(w,r))$. \\
$T_w$ & 9, 10 & For $w\in S$, the subgrid $\rho\Z \cap (Q(w,2r)\setminus Q(w,r))$. \\
$\Omega_1, \Omega_2, \Omega_3$ & 9, 10 & Global atypical events with negligible probability.\\
$\Omega_{4,w}$ & 9 & Local atypical event. \\
$\Omega_{5,W}$ & 9 & The event $\Omega_{4,w}$ will hold for more than half of $W$. \\
$\Omega_{6,W}$ & 9 & The event that more than  half of $W$ is percolating.\\
$\Omega_{7,W}$ & 9 & The event that for all $w\in W$ there are many black and not $4$-crowded points in $T_w$ and $\Omega_3^c\cap\Omega_{4,w}^c$ occurred.\\
$\pd(x,y)$ & 10 & The centered potential difference function. \\
$\pd(x,y|A)$ & 10 & Contribution to potential difference from stars in $A$.
\end{tabular}
\end{table}
\end{center}

\vspace{-30.0pt}
Lastly, the following lemma gathers some standard deviations
estimates on Poisson random variables; see \cite{kallenberg} for more details.

\begin{lem} \label{poissonlemma}
Let $X$ be a Poisson random variable with mean $\lambda$. Then:\\
(i) If $t\ge 2\lambda$ then $$\P(X\ge t)\le e^{-\frac{1}{4}t\log\left(\frac{t}{\lambda}\right)}. $$
(ii) There exists a $\delta>0$ such that for all $t\in[0,\delta]$ we have
$$ \P(|X-\lambda|\ge t \lambda) \le 2e^{-\lambda t^2/3}. $$
\end{lem}

\begin{proof}
For $t\ge \lambda$ set $s=t/\lambda$ in the inequality $s^t \P(X\ge t)\le \E(s^X)=e^{\lambda(s-1)}$, to
get
\begin{equation}\label{eq:poistar1}
\P(X\ge t) \le e^{-t \log\left(\frac{t}{\lambda}\right)+t-\lambda} = e^{-\lambda\left(\frac{t}{\lambda}
\log\left(\frac{t}{\lambda}\right)-\frac{t}{\lambda}+1\right)}.
\end{equation}
Since if $u:=\frac{t}{\lambda}\ge 2$, the inequality $\frac{3}{4}u\log u>u-1$ can be seen to hold, we get in that
case that $\P(X\ge t)\le e^{-\frac{1}{4}t\log\left(\frac{t}{\lambda}\right)}$, proving (i).

To prove (ii), note that by the same method, if $0<t\le \lambda$ one can set $s=t/\lambda$ in
the inequality $\P(X\le t)\le s^{-t}\E(s^X)$ to obtain
\begin{equation}\label{eq:poistar2}
\P(X\le t) \le e^{-\lambda\left(\frac{t}{\lambda}
\log\left(\frac{t}{\lambda}\right)-\frac{t}{\lambda}+1\right)}.
\end{equation}
Now, from \eqref{eq:poistar1} we get using a second order Taylor approximation that
$$ \P(X-\lambda \ge t\lambda) \le e^{-\lambda((1+t)\log(1+t)-t)}\le e^{-\lambda t^2/3} $$
for $t\in [0,\delta]$. A similar bound for $\P(X-\lambda\le t\lambda)$ follows similarly from \eqref{eq:poistar2}.
\end{proof}

\section{Proof outlines \label{sectionoutlines}}

We give a sketch of the proofs of the main results in the paper. This section is only included as an outline in order to give the reader a general feeling for the ideas used and to facilitate understanding of the detailed proofs in the later sections.

\bigskip
\noindent {\bf Equal volume of the basins of attraction.}
Of all the results mentioned above, one of the most interesting and surprising is that a.s.\  all the basins of attraction have volume $1$. This claim is relatively easy to prove, given the fact that the basins of attraction are a.s.\  bounded, and if we also assume that they have piecewise smooth boundaries (in Section \ref{sectionproofthm2} we give a detailed proof of the equal volumes property which does not use any information on the smoothness of the boundaries). Here is the proof, which is an adaptation of an argument due to Boris Tsirelson \cite{sodintsirelson}. First, we need an alternative expression for the force $F(x)$ that does not involve a different order of summation at every point $x$. In Section \ref{sectionproofprop1} we prove the following formula.

\begin{prop} \label{rearrangement}
Almost surely, for any $x\in\R^d\setminus {\cal Z}$ we have
\begin{equation}\label{eq:forcefield2}
F(x) = \sum_{z\in {\cal Z},\ |z|\uparrow}^\infty \frac{z-x}{|z-x|^d} + \kappa_d x,
\end{equation}
where $\kappa_d = \pi^{d/2}/\Gamma(d/2+1)$ is the volume of the unit ball in $\R^d$ and the summation is in order of increasing distance from $0$.
\end{prop}

Now, for a given basin of attraction $B(z_0)$, consider the oriented surface integral
$$ \int_{\partial B(z_0)} F(x)\cdot {\bf n}\ dS, $$
where ${\bf n} $ is the outward-pointing normal vector. We evaluate this integral in two ways. First, it is equal to $0$, since by the definition of the basin of attraction, on $\partial B(z_0)$ we actually have $F(x)\cdot {\bf n}=0$; this is because if $F(x)$ had a component in the direction of ${\bf n}$, there would be a flow curve crossing from one side of $\partial B(z_0)$ to the other. Second, the integral may be evaluated using the divergence theorem. Note that the function $g:\R^d\to \R^d$ defined by $g(x) = x/|x|^d$ satisfies $\textrm{div}(g) = d\kappa_d \delta_0$ in the distribution sense, where $\delta_y$ denotes a Dirac delta function at $y$. Therefore, using \eqref{eq:forcefield2} we have
$$ \textrm{div}(F) = -d\kappa_d \sum_{z\in{\cal Z}} \delta_z + d \kappa_d $$
in the distribution sense, and therefore, since $B(z_0)$ contains only the star $z_0$, we have
$$ 0 = \int_{\partial B(z_0)} F(x)\cdot {\bf n}\ dS = \int_{B(z_0)} \textrm{div}(F) dx = d\kappa_d (-1 + \vol(B(z_0))),$$
whence $\vol(B(z_0)) = 1$.

\paragraph{The proof of Theorem 3.} Our proof of Theorem 3 was inspired by, and follows the rough outline of, the proof of the main result of \cite[version 1]{nazarovetal2}, though several new conceptual and technical features are added. The basic idea is as follows. The event $E_R$ is defined in terms of the continuous-space force field $F(x)$ and is therefore hard to control. We bound it in terms of discrete events, by dividing space into a grid of cubes of side length $r \approx (\log R)^{2/d}$. Introduce a gravitational potential energy function $U(x)$ whose differences $U(x)-U(y)$ are the line integral of the gravitational force. Let $B=R^{8/9}$. If there is a gravitational flow curve $\Gamma$ crossing between $\partial Q(0,R)$ and $\partial Q(0,2R)$, then either $U(x)>B$ for some $x\in Q(0,2R)$ (an event which can be shown to be of negligible probability), or if not, then ``many'' (a positive fraction) of the $r$-grid cubes intersected by the curve $\Gamma$ have the property that either $U(x)\le -B$ for $x$ in that part of $\Gamma$ that intersects the cube or the change in potential energy along that part of the curve that intersects the cube is smaller than a constant times $Br/R$. Call such an $r$-cube ``bad''.

Now, if we could prove that the probability for each cube to be bad is bounded from above by $R^{-\delta}$ for some small $\delta>0$, and that the events of different cubes being bad are approximately independent on an appropriate scale, then Theorem \ref{thm3} would follow using standard subcritical percolation techniques. To bound the probability of a cube to be bad, we divide each cube into a grid of smaller cubes of size $\rho =R^{-1/10}$, and show that if the $r$-cube is bad, it contains many ``black'' $\rho$-subcubes, where, roughly, a subcube is called black if it contains a point $x$ where the norm of the force $F(x)$ is smaller than a constant times $B/R$. The probability of a cube to have many black subcubes is bounded using a first moment bound, which in this scale seems like the best one can do because of the extreme dependence of these events (since $\rho <<1$).

As for the approximate independence of the events of different cubes being bad, this is not strictly true in the scales under consideration. It is shown that the independence requirement can be replaced by a theorem bounding the joint density of the force field $F(x)$ evaluated at some set of points.

To make this skeleton of a proof work, several novel features are required. Detailed large deviations estimates are obtained for the gravitational potential, the force and its derivative. The joint distribution of the vector of values of forces at a given set of points is analyzed in detail. In the analysis of the bad cubes, it is necessary to bound the contributions to the potential energy
from two asymptotic regimes: First, from the effects of nearby stars causing the potential function to be close to $-\infty$; this is dealt with using a separate percolation argument. Second, from the ``intermediate'' range consisting of the scales between $r$ and $R^{1/d}$; these contributions are dealt with by diluting the set of potentially bad cubes by at most a factor $1/2$ and using a geometric covering lemma. Third, from the far range of stars at distance $> R^{1/d}$; this is dealt with using the large-deviation estimates.
In Dimensions 3 and 4, a more delicate argument is required involving a potential energy function that is not translation-invariant and has worse large-deviation behavior than in high dimensions.

\section{Derivation of Theorem \ref{thm2} \label{sectionproofthm2}}

We now show how Theorem \ref{thm2} follows from Theorem \ref{thm3} and Proposition \ref{rearrangement}.

First, the fact that $\p(E_R)\to 0$ as $R\to \infty$ clearly implies that a.s.\
all the basins of attraction are bounded.

Next, let $z_0 = \psi_{\cal Z}(0)$ be the star whose basin of attraction contains $0$, and let $B_0 = \psi_{\cal Z}^{-1}(\psi_{\cal Z}(0))$ be the basin of attraction of $z_0$. Let
$X$ be the allocation diameter $X=\text{diam}(B_0)$. If $X\ge R$ then
there exists an $x\in B_0$ with
$|x|\ge \frac{R}{2}$, and therefore $||x||_\infty \ge \frac{1}{2\sqrt{d}}R$.
Now, if $||z_0||_\infty\ge \frac{1}{4\sqrt{d}}R$ then since $0$ is in the basin of attraction of $z_0$ the event $E_{R/8\sqrt{d}}$ happened. Otherwise, since $x$ is in the basin of attraction of $z_0$ and we have $||z_0||_\infty < \frac{1}{4\sqrt{d}}R < \frac{1}{2\sqrt{d}} R \le ||x||_\infty$, the event $E_{R/4\sqrt{d}}$ happened. So we have shown that
$$ \p(X\ge R) \le \p(E_{R/8\sqrt{d}}) + \p(E_{R/4\sqrt{d}}). $$
This implies that the estimates \eqref{eq:tailestimate1} and \eqref{eq:tailestimate2} follow from the corresponding estimates in Theorem \ref{thm3}.

Next, we show that a.s.\ all basins of attraction have volume 1. We
use a variant of the argument sketched in Section
\ref{sectionoutlines}
which does not require any knowledge about the smoothness of the boundary of $B_i$. A similar argument in a slightly different context was briefly mentioned in \cite[Section 12.2]{nazarovetal2}.

Let $B_i$ be the basin of attraction of the star $z_i$.
As in the introduction, for $x\in B_i$ denote by $\tau_x$ the time that it takes $x$ to flow into $z_i$, or equivalently the maximum time up to which the integral curve $Y_x(\cdot)$ is defined; for continuity we set $\tau_{z_i} = 0$.
For $0\le a\le b\le \infty$ denote
\begin{eqnarray*}
E_{a,b} &=& \{ x\in B_i\,:\,a\le\tau_x\le b \}, \\
V_{a,b} &=& \vol(E_{a,b}).
\end{eqnarray*}
Note that on $B_i\setminus\{z_i\}$ the force field $F$ satisfies
$\text{div}(F) \equiv -d\kappa_d$, by taking the divergence of each term
in \eqref{eq:forcefield2}
(see Lemma \ref{simconv1} below, which justifies interchanging the divergence and summation operations). Therefore,
by a version of Liouville's theorem \cite[Lemma 1, p. 69]{arnold}, it follows
that $\frac{d}{dt}V_{t,\infty} = -d\kappa_d V_{t,\infty}$, so
$$ V_{t,\infty} = V_{0,\infty}e^{-d \kappa_d t}. $$
In particular, for $t\searrow 0$ we get that
\begin{equation}\label{eq:equating1}
V_{0,t} = V_{0,\infty}-V_{t,\infty} = V_{0,\infty}\left(1-e^{-d\kappa_d t}\right)
= d\kappa_d V_{0,\infty}t + O(t^2).
\end{equation}
Estimate $V_{0,t}$ for small $t$
in a different way, as follows. In a neighborhood of $z_i$ the field
$F$ satisfies $F(x) = \frac{z_i-x}{|z_i-x|^d} + O(1)$. Without the error term
it would be easy to solve the differential equation explicitly, so this
implies by approximation that
$$ B(z_i,(dt)^{1/d}-o(t^{1/d})) \subseteq E_{0,t} \subseteq B(z_i,(dt)^{1/d} + o(t^{1/d})) $$
(this relies on the following easily-verified claim regarding a
one-dimensional differential inequality: if $g$ is a real-valued function
on $[0,\infty)$ that satisfies $g(0)=0, |g'(t)-g(t)^{1-d}|\le C$, then
$g(t)=(dt)^{1/d}+o(t^{1/d})$ when $t\searrow 0$),
and therefore that
\begin{equation}\label{eq:equating2}
V_{0,t} = d\kappa_d t + o(t).
\end{equation}
Equating \eqref{eq:equating1} and \eqref{eq:equating2} gives that
$ \vol(B_i) = V_{0,\infty} = 1$, as claimed.

We have shown that the basins of attraction are a.s. all bounded and have volume 1, and they are clearly disjoint. The last claim that needs proving is that a.s. they cover all of $\R^d$ except a set of
measure $0$. We use the following {\bf mass transport lemma}.

\begin{lem} \label{transport}
Let $f:\Z^d\times \Z^d\to[0,\infty)$ satisfy $f(m,n) = f(m+u,n+u)$ for any $m,n,u\in \Z^d$. Then
for all $n\in\Z^d$ we have that
$$ \sum_{m\in\Z^d} f(m,n) = \sum_{m\in\Z^d} f(n,m). $$
\end{lem}

\begin{proof}
$f(m,n)=f(m-n,0)=:g(m-n)$, and both sums become just $\sum_n g(n)$.
\end{proof}

Define $f:\Z^d\times \Z^d\to\R$ by
$$ f(m,n) = \E\bigg[\vol\bigg(Q(m,1/2) \cap \bigcup_{i\,:\,z_i\in Q(n,1/2)} B_i\bigg)\bigg]. $$
or in words the expected volume of the part of $Q(m,1/2)$ that gets allocated to some $z_i \in Q(n,1/2)$.
Note that $\sum_{m\in \Z^d} f(m,n)$ represents the expected volume of points in $\R^d$ being
allocated to some $z_i \in Q(n,1/2)$. Since we showed that $\vol(B_i)=1$ for all $i$, this is
equal to the expected number of $z_i\in Q(n,1/2)$, which is $\vol(Q(n,1/2))=1$. So, if we denote
$D = \cup_{i=1}^\infty B_i$, the union of all the basins of attraction, then
 by Lemma \ref{transport} we get that for all $n\in\Z^d$ we have
$$ 1 = \sum_{m\in \Z^d} f(n,m) = \E\Big[\vol(D\cap Q(n,1/2))\Big]. $$
The random variable $\vol(D\cap Q(n,1/2))$ is bounded from above by 1. If its expected value is 1 then it is 1 almost surely. Therefore almost surely
$$ \vol(\R^d\setminus D) = \sum_{n\in \Z^d} \vol(Q(n,1/2)\setminus D) = \sum_{n\in Z^d} 0 = 0, $$
as claimed.
\qed

\section{Existence and differentiability of $F$ \label{sectionproofprop1}}

In this section we prove Propositions \ref{prop1} and
\ref{rearrangement}.

\subsection{Proof of a.s.\  convergence of $F(x)$ \label{subsectionasconv}}

First, let us prove that $F(x)$ is well-defined, that is, that the sum in \eqref{eq:forcefield} converges a.s.\  for fixed $x\in\R^d$. Since the sum is defined in a translation-invariant manner, it is clearly enough to prove that the sum for $F(0)$ converges a.s. Let $\rho_0=0$, and for $i\ge 1$ let $\rho_i = |z_i|$ be the distance of $z_i$ from the origin. Since ${\cal Z}=(z_i)_i$ is a Poisson process, we then have that the random variables $\big(\kappa_d(\rho_i^d-\rho_{i-1}^d)\big)_{i=1}^\infty$ are i.i.d.\ with $\text{Exp}(1)$ distribution (recall $\kappa_d=\vol(B(0,1))$), and therefore by the law of large numbers, almost surely
$$ \ \ \ \ \qquad\qquad
\frac{\rho_i}{i^{1/d}} = \left[ \frac{\sum_{j=1}^i (\rho_j^d-\rho_{j-1}^d)}{i} \right]^{1/d}
 \xrightarrow[i\to\infty]{} \kappa_d^{-1/d}. \qquad \qquad \text{(LLN)} $$

Now, if we condition on the values of $(\rho_i)_i$, thinking of them as a deterministic sequence such that $\rho_i/i^{1/d} \to \kappa_d^{-1/d}$ as $i\to\infty$, then each $z_i$ is distributed uniformly on the sphere of radius $\rho_i$ around the origin. For any $i\ge 1$, each term $(z_i)/|z_i|^d$ in the sum in \eqref{eq:forcefield} (where $x$ is taken as $0$) has (conditional) mean 0 and variance bounded by
$O\big(\rho_i^{-2(d-1)}\big) = O\left(i^{-2(d-1)/d}\right)$. Since in the event (LLN) the sum of the variances converges (note that this fails in dimension 2), by the Kolmogorov three-series theorem the sum in \eqref{eq:forcefield} converges a.s. This is true a.s.\  conditionally on $(\rho_i)_i$, therefore it is true a.s.\  and $F(0)$ is defined.

\subsection{Simultaneous convergence with a fixed order of summands
\label{subsectionproofsimult}}

Denote as before $\rho_i=|z_i|$.
Let $g:\R^d\to\R^d$ be defined by $g(x)=x/|x|^d$.
Let $i_0 = \min\{i\,:\,\rho_i\ge 2\}$. Define
\begin{equation}\label{eq:sumhx}
H(x) = \sum_{i=i_0}^\infty
\frac{z_i-x}{|z_i-x|^d} = \sum_{i=i_0}^\infty g(z_i-x),
\end{equation}
a function that we will see shortly is closely related to $F(x)$.

\begin{lem} \label{simconv1}
Almost surely, the sum defining $H(x)$ converges simultaneously
and uniformly for all $x\in B(0,1)$
and defines a continuously differentiable function. The series can be differentiated
termwise.
\end{lem}

\begin{proof}
For a function $f:\R^i\to\R^j$, denote by $D_k f$ the tensor of $k$-th derivatives
of $f$ (which can be thought of as a $(j\cdot i^k)$-dimensional vector). Note
that for $|x|>1$ we have that
\begin{eqnarray}
| D_1 g(x) | &=& O( |x|^{-d} ), \label{eq:apriori1} \\
| D_2 g(x) | &=& O( |x|^{-d-1} ), \label{eq:apriori2}
\end{eqnarray}
and in general for any $k\ge0$ we have
\begin{equation} \label{eq:apriori3}
 | D_k g(x) | = O(|x|^{-d-k+1} ),
\end{equation}
where the constant implicit in the big-$O$ depends on $d$ and on $k$.
The best way to see \eqref{eq:apriori1} is to represent $D_1 g(x)$, the matrix
of the first differential of $g$ at $x$, in an orthonormal coordinate system
containing the radial direction $x/|x|$; this gives a diagonal matrix whose
entries are $d-1$ copies of $|x|^{-d}$ and one copy of
$-(d-1)|x|^{-d} = \frac{d}{dr}\big|_{r=|x|} r^{-(d-1)}$, so in fact
$|D_1 g(x)| = C_{1}|x|^{-d}$, where $C_{1}=\big(d(d-1)\big)^{1/2}$.
Equations \eqref{eq:apriori2} and \eqref{eq:apriori3} can be proved similarly.

Now, similarly to \eqref{eq:sumhx}, define for $x\in B(0,1)$
\begin{eqnarray}\label{eq:sumhx1}
H_1(x) &=& \sum_{i=i_0}^\infty D_1[g(z_i-x)], \\
H_2(x) &=& \sum_{i=i_0}^\infty D_2[g(z_i-x)]. \label{eq:sumhx2}
\end{eqnarray}
Condition on the $\rho_i$, and condition on the event (LLN). In the previous
subsection we showed that the sum defining $H(0)$ converges a.s., so assume that this holds.
A similar argument shows that the sum defining $H_1(0)$ converges a.s., so
condition on that as well. We shall show that under these conditions
the sum in \eqref{eq:sumhx} converges uniformly on $B(0,1)$ to a $C^1$ function.

First, \eqref{eq:apriori2} together with the assumption that (LLN) holds
immediately imply that the sum \eqref{eq:sumhx2}
converges absolutely uniformly on $B(0,1)$, and similarly from \eqref{eq:apriori3}
the same is true for the sums of differentials of all orders $k\ge 2$. In
particular it follows that (under the above conditioning)
$H_2(x)$ is a $C^\infty$ function on $B(0,1)$.

Next, $H_1(x)-H_1(0)$ can be represented as a sum of line integrals from $0$ to $x$
of the terms in the sum for $H_2(x)$. Therefore the sum for $H_1(x)-H_1(0)$
converges uniformly on $B(0,1)$ to a function whose differential is $H_2(x)$,
and since we assumed that the sum for $H_1(0)$ converges, it follows that
the sum \eqref{eq:sumhx1} converges uniformly on $B(0,1)$ to a differentiable
function. Similarly $H(x)-H(0)$ can be represented as a sum of line integrals
of the terms in \eqref{eq:sumhx1}, so repeating the above argument, using the fact
that we assumed that the sum for $H(0)$ converges, gives that the sum in \eqref{eq:sumhx}
converges uniformly on $B(0,1)$ to a $C^1$ (in fact, $C^\infty$) function. This was
true under the conditioning on an almost sure event, so the lemma is proved.
\end{proof}

\subsection{The rearrangement identity \label{subsectionrearrange}}

If $u,x\in\R^d$ we denote
$$ G^{\{u\}}(x) = \sum_{|z_i-u| \uparrow} \frac{z_i-x}{|z_i-x|^d} $$
(the terms are summed in order of increasing distance from $u$, and
this sum is defined a.s.\ as with $F(x)$).
\begin{lem} \label{rearr}
For any $x,u,v\in\R^d$ we have that a.s.
\begin{equation}\label{eq:rearrangeidentity}
G^{\{u\}}(x)-G^{\{v\}}(x) = \kappa_d(u-v).
\end{equation}
\end{lem}

\begin{figure}
\centering \resizebox{!}{6cm}{\includegraphics{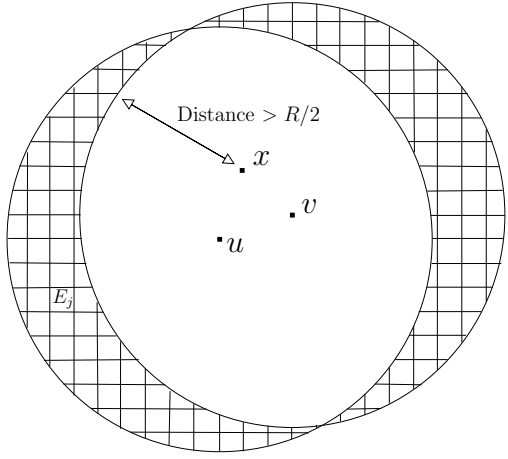}}

\vspace{-140.0pt} \hspace{-10.0pt} \resizebox{!}{0.21 cm}{Distance $>R/2$}

\vspace{2.0pt} \hspace{0.0pt} {$x$}

\vspace{4.0pt} \hspace{36.0pt} {$v$}

\vspace{0.0pt} \hspace{-20.0pt} {$u$}

\vspace{5.0pt} \hspace{-150.0pt} \resizebox{!}{0.18 cm}{$E_j$}

\vspace{60.0pt}

\caption{The balls $B(u,R)$ and $B(v,R)$ and the sets $E_j$. \label{figrearrange}}
\end{figure}

\begin{proof} First, compute expectations: Let $N_{u,x}$ be the (random) number of stars in the ball $B(u,|u-x|)$. Recall the well-known physics principle that says that the total gravitational pull on a point $x$ from a uniformly distributed spherical shell of mass with center $u$, radius $r$ and total mass $M$ is equal to $0$ if $r> |u-x|$ and to $M(u-x)/|u-x|^d$ if $r < |u-x|$ (this last fact follows from the harmonicity of the function $x\to(u-x)/|u-x|^d$). Therefore, by conditioning on the distances of the stars from $u$ (as was done in Section \ref{subsectionasconv} above with $u=0$), we get that
$$ \E\Big[ G^{\{u\}}(x) \Big| N_{u,x} \Big] = N_{u,x}\cdot \frac{u-x}{|u-x|^d}. $$
Therefore
$$ \E\Big[ G^{\{u\}}(x) \Big] = \E[N_{u,x}] \frac{u-x}{|u-x|^d} = \kappa_d|u-x|^d \frac{u-x}{|u-x|^d}
= \kappa_d(u-x), $$ so
$$ \E\Big[ G^{\{u\}}(x) \Big] - \E\Big[ G^{\{v\}}(x) \Big] = \kappa_d(u-v). $$
Now, let $R>0$ be large, and consider the truncated series
$$ G_R^{\{u\}}(x) = \sum_{|z_i-u|<R} \frac{z_i-x}{|z_i-x|^d}. $$
Then $$ G_R^{\{u\}}(x) - G_R^{\{v\}}(x) = \sum_{z_i \in A_R} \frac{z_i-x}{|z_i-x|^d}
- \sum_{z_i \in B_R} \frac{z_i-x}{|z_i-x|^d}, $$
where $A_R = B(u,R) \setminus B(v,R), B_R = B(v,R) \setminus B(u,R). $ We show that the variance
of this expression tends to $0$ when $R\to\infty$: Partition the set $B(u,R)\triangle B(v,R)$ into
$O(R^{d-1})$ disjoint sets $(E_j)_j$ of Lebesgue measure $O(1)$ such that each $E_j$ is contained in either $A_R$ or $B_R$ (see Figure \ref{figrearrange}; the constant in the big-$O$ depends on $u$ and $v$), and for each $j$ let
$$ Y_j = \sum_{z_i\in E_j} \frac{z_i-x}{|z_i-x|^d} $$ be the contribution to the force from stars in $E_j$.
Then we can write
$$ G_R^{\{u\}}(x) - G_R^{\{v\}}(x) = \sum_{j} \pm Y_j. $$
The $Y_j$'s are independent, and each has variance bounded from above by
\begin{eqnarray*}
\E|Y_j|^2 &=& \E\Big[ \E\big[ |Y_j|^2 \big| \textrm{card}(E_j \cap \{z_i\}_i) \big] \Big] \\ &
\le & \frac{C}{R^{2d-2}} \E\Big[ \textrm{card}(E_j \cap \{z_i\}_i) \Big] = O\left( \frac{1}{R^{2d-2}} \right).
\end{eqnarray*}
(Note that this is true since $B(u,R)\triangle B(v,R) \subset \R^d \setminus B(x,R/2)$ for sufficiently large $R$, see Figure \ref{figrearrange}). Therefore
$$ \textrm{Var}\left( G_R^{\{u\}}(x) - G_R^{\{v\}}(x) \right) = O(R^{-d+1}) \xrightarrow[R\to\infty]{} 0, $$
which finishes the proof, since a.s.\  $G_R^{\{u\}}(x)\to G^{\{u\}}(x)$ and $G_R^{\{v\}}(x)\to G^{\{v\}}(x)$ as $R\to\infty$
\end{proof}

\subsection{Proof of Propositions \ref{prop1} and \ref{rearrangement}}

Both Propositions \ref{prop1} and \ref{rearrangement} follow
immediately from the following theorem.

\begin{thm} \label{thmrearrange}
With probability 1, the following four statements hold: \\
(i) The sum defining $G^{\{u\}}(x)$ converges simultaneously for all $u\in \R^d$
and $x\in\R^d\setminus{\cal Z}$. \\
(ii) The convergence is uniform on compact sets in $\R^d\times(\R^d\setminus {\cal Z})$.
\\
(iii) The rearrangement identity
\eqref{eq:rearrangeidentity} holds for all $u, v$ and $x$. \\
(iv) For all $u$ the function $G^{\{u\}}(x)$ is continuously differentiable in $x$.
\end{thm}

\begin{proof}
By Lemma \ref{simconv1} we know that, off of a null event $\Omega_1$, for all rational
$u\in\R^d$ the sum defining $G^{\{u\}}(x)$ converges simultaneously
for all $x\in\R^d\setminus{\cal Z}$, and the convergence
is uniform for $x$ ranging in a compact set.

By Lemma \ref{rearr}, off of a null event $\Omega_2$, the identity
\eqref{eq:rearrangeidentity} holds
provided $u,v,x$ are rational. (By the continuity in $x$, the assumption that $x$
is rational can be dropped outside $\Omega_1\cup\Omega_2$.)

Let $N(u,R,\eps)$ be the number of stars in the $\eps$-neighborhood of the sphere
of radius $R$ around $u$. The mean of $N(u,R,\eps)$ is at most $C_d R^{d-1} \eps$.
Let $A(u,R,\eps)$ be the event that $N(u,R,\eps)$ is less than twice its mean. Then
by Lemma \ref{poissonlemma},
\begin{equation}\label{eq:upbound}
\p( A(u,R,\eps)^c ) \le \exp(-a R^{d-1}\eps)
\end{equation} for some $a>0$. It follows by
Borel-Cantelli that off of a null event $\Omega_3$, for each rational $q$ and $\eps>0$,
there is a (random) $R_* = R_*(q,\eps)$ such that $A(q,R,\eps)$ holds for all
$R>R_*(q,\eps)$ that are multiples of $\eps$.

Now fix a configuration of stars $\omega \notin \Omega_1\cup\Omega_2
\cup\Omega_3$, and choose $\eps>0$. For each $u\in\R^d$ find a rational $q=q(u,\eps)$
within distance $\eps$ of $u$. Then for $R>R_*(q,\eps)$ and $x\in B(u,R/3)$,
we have
$$ \left|G_R^{\{u\}}(x) - G_R^{\{q\}}(x)\right| < 2N(u,R,\eps)(2/R)^{d-1}
<2^{d+1}C_d \eps $$
where $G_R^{\{u\}}(x)$ is defined as in the proof of Lemma \ref{rearr} above.
Thus
$$ \limsup_R \left|G_R^{\{u\}}(x)- G_R^{\{q\}}(x)\right| \le
2^{d+1} C_d \eps$$
for all $\eps>0$. This verifies (i) for $\omega$, and (ii), (iii), (iv)
follow similarly by approximation.
\end{proof}

\section{The joint density of a vector of forces \label{sectionforcethm}}

In this section we prove an estimate that will be required in the proof
of Theorem \ref{thm3}. Suppose we have $N$
points $x_1,\ldots,x_N\in\R^d$ with $|x_i-x_j|> S$ for every $i\ne j$.
Fix a positive $\lambda$, and define the event
\begin{equation*}
E=\Big\{\text{There is at least one star in }B(x_i,\lambda)\text{ for every }
1\le i\le N\Big\}.
\end{equation*}
Denote by ${\cal M}$ the $\sigma$-algebra generated\ \ by the locations
of the stars in
$\left(\cup_{i=1}^N B(x_i,\lambda)\right)^c$. Denote by $X$ the random vector of
forces $(F(x_i))_{1\le i\le N}$.
Then we have the following bound on the joint density of $X$.
\begin{thm} \label{joint_force_density_thm}
There exist constants $c_0, C_1>0$ (depending on the dimension $d$) such that
if
\begin{equation}\label{eq:assumption}
\lambda<c_0\frac{S}{(\log N)^{\frac{1}{d}}},
\end{equation}
then conditioned on the
event $E$ and on the $\sigma$-algebra ${\cal M}$, almost surely
the joint density of $X$ exists and is bounded from above by
$ (C_1 \lambda^{d^2-d})^N$.
\end{thm}

We will use the following two simple lemmas.

\begin{lem} \label{simple1} There exists a constant $C_7>0$ (depending on $d$) such that
if $x_1,x_2,\ldots,x_N\in \R^d$ satisfy $|x_i-x_j|>S$ for all $i\neq j$, then
for all $1\le i\le N$ we have that
$$\sum_{j=2}^N \frac{1}{|x_j-x_1|^d} \le \frac{C_7 \log N}{S^d}. $$
\end{lem}

\begin{proof}
For $k=1,2,\ldots$, let $V_k = \{x_i\}_{i=2}^N \cap \Big(B(x_1,2^k S)\setminus B(x_1,2^{k-1} S)\Big)$.
Clearly $|V_k| \le \frac{\vol(B(x_1,(2^k+1)S))}{\vol(B(0,S/2))} \le 4^d\cdot 2^{kd}$, and also trivially
$|V_k| \le N$. Therefore we have
\begin{eqnarray*}
\sum_{j=2}^N \frac{1}{|x_j-x_1|^d} &=& \sum_{k=1}^\infty \sum_{x\in V_k} \frac{1}{|x-x_1|^d}
\le \sum_{k=1}^\infty \frac{\min(N, 4^d 2^{kd})}{2^{kd} S^d/2^d} \\
&=& 8^d \sum_{k=1}^{\lfloor \frac{\log N}{d \log 2} \rfloor-1} \frac{1}{S^d} + 2^d N
\sum_{k=\lfloor \frac{\log N}{d \log 2} \rfloor}^\infty \frac{1}{2^{kd} S^d} \\ &=& O\left(\frac{\log N}{S^d}\right)
+ O\left(\frac{1}{S^d}\right) = O\left(\frac{\log N}{S^d}\right).
\end{eqnarray*}
\end{proof}

\begin{lem} \label{simple2}
If $A=(a_{i,j})_{i,j=1}^k$ is a matrix such that
$|a_{i,i}|\ge 2 \sum_{j\neq i} |a_{i,j}|$ for all $1\le i\le k$, then
$$ |\det A| \ge \frac{1}{2^k} \prod_{i=1}^k |a_{i,i}|. $$
\end{lem}

\begin{proof} This is a variant of Hadamard's theorem in linear algebra.
First, by multiplying each row of $A$ by $a_{i,i}^{-1}$, we may assume without
loss of generality that $a_{i,i}=1$ for all $1\le i\le k$, so $A=I+B$,
where $I$ is the identity matrix and $B=(b_{i,j})$ is a matrix such that
$\sum_j |b_{i,j}| \le 1/2$ for all $1\le i\le k$. For all $x\in \R^k$
we have
$$||Bx||_\infty \le ||x||_\infty \max_i \sum_j |b_{i,j}|
\le \frac{1}{2}||x||_\infty. $$
Therefore $$||Ax||_\infty = ||Ix+Bx||_\infty \ge ||x||_\infty - ||Bx||_\infty
\ge \frac{1}{2} ||x||_\infty.$$
We have shown that all the eigenvalues of $A$ are greater in absolute value
than $1/2$, therefore $|\det A| \ge 2^{-k}$, as claimed.
\end{proof}

\begin{proof}[Proof of Theorem \ref{joint_force_density_thm}]
Let us condition everything on the event $E$ and moreover on the number
of stars $\nu_i$
in $B(x_i,\lambda)$ for each $1\le i\le N$. With this conditioning,
the set of stars in $B(x_i,\lambda)$ is simply a vector
$(Y_{i,1},Y_{i,2},\ldots,Y_{i,\nu_i})$ of $\nu_i$ i.i.d.\
points chosen uniformly in $B(x_i,\lambda)$. Now condition further on
the $\sigma$-algebra ${\cal M}$ and on the $\sigma$-algebra
${\cal L}$ generated by the locations of the stars
\hbox{$\big\{ Y_{i,j}\ |\ 1\le i\le N, 2\le j\le \nu_i\big\}$.}
This leaves only the stars $(Y_{i,1})_{1\le i\le N}$ as a source of randomness.
The vector $X$ of forces can therefore be written as
$$ X = G(Y_{1,1},Y_{2,1},\ldots,Y_{N,1}) + Z, $$
where $Z$ represents the contribution that is measurable with respect
to ${\cal M}\vee{\cal L}$, and where
$G:B(x_1,\lambda)\times B(x_2,\lambda)\times \ldots \times B(x_N,\lambda)
\to\R^{Nd}$ is the function defined by
$$ G(y_1,\ldots,y_N) = \left(\ \sum_{j=1}^N \frac{y_j-x_1}{|y_j-x_1|^d},\
\sum_{j=1}^N \frac{y_j-x_2}{|y_j-x_2|^d},\
\ldots\ ,\ \sum_{j=1}^N \frac{y_j-x_N}{|y_j-x_N|^d}
                       \right).
$$
Denote ${\bf B}=
B(x_1,\lambda)\times B(x_2,\lambda)\times \ldots \times B(x_N,\lambda)$.
The volume of ${\bf B}$ is $\kappa_d^N \lambda^{Nd}$.
Therefore, to prove that the joint density of $G(Y_{1,1},Y_{2,1},\ldots,Y_{N,1})$, and therefore
also the joint density of $X$ conditioned on the event $E$ and on the $\sigma$-algebra
${\cal M}\vee {\cal L}$, is bounded from above by $(C_1 \lambda^{d^2-d})^N$,
it will be enough to prove two things: First, that the function $G:{\bf B}\to \R^{Nd}$ is one-to-one; and second, that the Jacobian of the function $G:{\bf B}\to \R^{Nd}$ is bounded from below by
$(C_2 \lambda^{d^2})^{-N}$, where $C_2>0$ is some large constant. Interestingly, both of
these claims require that the assumption \eqref{eq:assumption} hold for some constant $c_0>0$.

We prove the first claim. Assume \eqref{eq:assumption}, where $c_0>0$ is some small constant whose value will be specified soon. Denote as before $g(x)=x/|x|^d$.
Fix $a_1,\ldots,a_N\in \R^d$. Our goal is to prove that if the system of equations
\begin{equation}\label{eq:system}
\sum_{j=1}^N g(y_j-x_i) = a_i, \quad 1\le i\le N
\end{equation}
has a solution $(y_1,\ldots,y_N) \in {\bf B}$, then this solution is unique.
The following proof of this fact was suggested by the referee, and simplified an
earlier proof.
Assume the contrary: $y=(y_1,\ldots,y_N)\in {\bf B}$, $y' = (y_1',\ldots,y_N')\in {\bf B}$, and
\begin{equation}\label{eq:appen1}
\sum_{j=1}^N g(y_j-x_i) = \sum_{j=1}^N g(y_j'-x_i)\ \ \text{ for }i=1,\ldots,N,
\end{equation}
and the number $\varepsilon = \max_j |y_j-y_j'|$ does not vanish. Without loss of generality assume that $|y_1-y_1'|=\varepsilon$. Introducing $u=g(y_1-x_1), u'=g(y_1'-x_1), v=\sum_{j=2}^N g(y_j-x_1), v'=\sum_{j=2}^N g(y_j'-x_1)$, we have $u+v=u'+v'$ (by \eqref{eq:appen1} for $i=1$; other $i$ will not be used), therefore
\begin{equation}\label{eq:appen2} |u-u'|=|v-v'|. \end{equation}
We note that $|y_j-x_1|\ge |x_j-x_1|-|y_j-x_j|\ge |x_j-x_1|-\lambda$ for $j>1$. Taking into account that $|x_j-x_1|\ge S\ge 2\lambda$ by \eqref{eq:assumption}, we get $|y_j-x_1|\ge \frac{1}{2}|x_j-x_1|$. By \eqref{eq:apriori1}, for $j>1$,
\begin{eqnarray*}
|g(y_j-x_1)-g(y_j'-x_1)|&\le& C \frac{|y_j-y_j'|}{\min(|y_j-x_1|^d,|y_j'-x_1|^d)}
\\ &\le& C\frac{|y_j-y_j'|}{|x_j-x_1|^d}\le C\frac{\varepsilon}{|x_j-x_1|^d}.
\end{eqnarray*}
By Lemma \ref{simple1},
$$ |v-v'|\le\sum_{j=2}^N |g(y_j-x_1)-g(y_j'-x_1)|\le C\varepsilon\sum_{j=2}^N\frac{1}{|x_j-x_1|^d}
\le \frac{C_{99} \varepsilon\log N}{S^d}. $$
Here $C_{99}$ does not depend on $c_0$ in \eqref{eq:assumption} as long as $\frac{c_0}{(\log N)^{1/d}}\le 0.5. $

Using \eqref{eq:appen2} and \eqref{eq:assumption}, we get that $|u-u'|=|v-v'|\le \frac{C_{99}\log N}{S^d}
\le C_{99}\varepsilon(c_0/\lambda)^d$, that is,
\begin{equation}\label{eq:appen3}
\frac{|g(y_1-x_1)-g(y_1'-x_1)|}{|y_1-y_1'|}\le C_{99}\left(\frac{c_0}{\lambda}\right)^d.
\end{equation}
The function $g:\R^d\setminus\{0\}\to \R^d\setminus\{0\}$ is invertible; $g^{-1}(z)=z/|z|^{\frac{d}{d-1}}$. Similarly to \eqref{eq:apriori1}, $|D_1 g^{-1}(z)|=O(|z|^{-\frac{d}{d-1}})$. Thus if we restrict $g$ to $B(0,\lambda)$ and accordingly $g^{-1}$ to $\R^d\setminus B(0,\frac{1}{\lambda^{d-1}})$, then $g^{-1}$ satisfies the Lipschitz condition with the constant $C_{98}\lambda^d$ (even though $\R^d\setminus B(0,\frac{1}{\lambda^{d-1}})$ is not convex...). We consider $z=g(y_1-x_1), z'=g(y_1'-x_1)$ and get by \eqref{eq:appen3} and the Lipschitz condition
$$ \frac{|z-z'|}{|g^{-1}(z)-g^{-1}(z')|}\le C_{99}\left(\frac{c_0}{\lambda}\right)^d, \qquad
\frac{|g^{-1}(z)-g^{-1}(z')|}{|z-z'|}\le C_{98}\lambda^d, $$
which is a contradiction if $c_0$ is small enough, namely if $C_{98}C_{99}c_0^d<1$.
This finishes the proof that $G$ is one-to-one.

It remains to prove that the Jacobian of $G$ is bounded from below by
$(C_2 \lambda^{d^2})^{-N}$, for some large constant $C_2>0$,
again assuming \eqref{eq:assumption}. The Jacobian matrix $J$
of $G$ can be written as a block matrix $(J_{i,j})_{1\le i,j\le N}$, where
each $J_{i,j}$ is the $d\times d$ Jacobian matrix of the function
$ y_j \to \frac{y_j-x_i}{|y_j-x_i|^d}$. Again by the computation
of the matrix $D_1 g(x)$, each $J_{i,j}$ is a diagonalizable
matrix with one eigenvalue equal to $-(d-1)|y_j-x_i|^{-d}$ and $d-1$
eigenvalues equal to $|y_j-x_i|^{-d}$. Furthermore, by choosing for each $y_i$
the appropriate radial coordinate system (as a function of $y_i$), we may
assume that the blocks $J_{i,i}$, $1\le i\le N$ are in diagonal form.
Any other block $J_{i,j}=(a_{i,j,k,l})_{1\le k,l\le d}$
for $i\neq j$ is not necessarily in diagonal
form, but its entries satisfy
$$ a_{i,j,k,l} \le C_5 |y_j-x_i|^{-d} \le C_6 |x_j-x_i|^{-d}. $$
Recall that our assumptions are that $|y_i-x_i| < \lambda$ and
$|y_j-x_i|>S-\lambda > S/2$ for $i\neq j$. We wish to apply Lemma
\ref{simple2} to the matrix $J$. By Lemma \ref{simple1}, the
assumptions of Lemma \ref{simple2} will hold if we have
$$\lambda^{-d} > 2 d C_6 C_7 \frac{\log N}{S^d}. $$
This holds if $\lambda < C_8 \frac{S}{(\log N)^{1/d}}$,
where $C_8 = (2d C_6 C_7)^{-1/d}$. The conclusion of Lemma \ref{simple2}
is exactly our claim.
\end{proof}

\section{The gravitational potential function \label{sectionpotential}}

We define a new function, the {\bf gravitational potential function}. It will be defined in
dimensions 5 and higher only, and is designed to be a stationary centered random function that has as minus its gradient the force function $F$. If $d\ge 5$, the gravitational potential function $U:\R^d\to \R$ is defined by
\begin{equation}\label{eq:potential}
U(x) = \frac{1}{d-2}
\lim_{T\to\infty} \Big[\sum_{i\,:\,|z_i-x|<T} \frac{-1}{|z_i-x|^{d-2}} + \frac{d\kappa_d}{2}T^2\Big].
\end{equation}


As with the case of the force, we need to check that the potential function is a.s. defined and is
well-behaved (in fact, in the case of the potential this is only true in dimensions 5 and higher).
For any $p > q \ge 0$, denote by $N_{p,q}$ the random number of stars in $B(0,p)\setminus B(0,q)$, and denote
$$U_{p,q} = \sum_{i\ :\ q<|z_i|\le p} \frac{1}{|z_i|^{d-2}}. $$
Let $W_{p,q}$ be a random vector distributed uniformly on $B(0,p)\setminus B(0,q)$.  An easy computation gives the following.
\begin{eqnarray}
\E\left[ |W_{p,q}|^\alpha \right] &=& \frac{d}{d+\alpha}\frac{p^{d+\alpha}-q^{d+\alpha}}{p^d-q^d},\qquad
  (\alpha\neq -d), \nonumber \\
\E[N_{p,q}] &=& \var[N_{p,q}] = \kappa_d(p^d-q^d), \nonumber \\
\E[U_{p,q}\ |\ N_{p,q}] &=& N_{p,q}\E\left[|W_{p,q}|^{2-d}\right] = N_{p,q}\cdot \frac{d}{2}\cdot
  \frac{p^2-q^2}{p^d-q^d}, \nonumber \\
\E[U_{p,q}] &=& \frac{d\kappa_d}{2}(p^2-q^2), \label{eq:expectationupq} \\
\var[U_{p,q}\ |\ N_{p,q}] &=& N_{p,q}\cdot \var\left[ |W_{p,q}|^{2-d}\right]  \nonumber \\ &=&
  N_{p,q}\cdot \left(\frac{d}{d-4}\cdot\frac{q^{4-d}-p^{4-d}}{p^d-q^d}-
    \frac{d^2}{4}\cdot\left( \frac{p^2-q^2}{p^d-q^d}\right)^2 \right), \nonumber \\
\var[U_{p,q}] &=& \E \Big[\var[U_{p,q}\ |\ N_{p,q}]\Big] + \var \Big[\E [U_{p,q}\ |\ N_{p,q}]\Big] \nonumber \\
 &=& \frac{d \kappa_d}{d-4}\left( \frac{1}{q^{d-4}} - \frac{1}{p^{d-4}} \right). \label{eq:varianceupq}
\end{eqnarray}
Now, from \eqref{eq:expectationupq} and \eqref{eq:varianceupq} it immediately
follows that when $d\ge 5$, for any fixed $x\in \R^d$ the limit in \eqref{eq:potential}
exists a.s.\  and defines a centered random variable.

For any $u,x\in\R^d$ denote
$$ H_R^{\{u\}}(x) = \frac{1}{d-2}
 \Big[\sum_{i\,:\,|z_i-u|<R} \frac{-1}{|z_i-x|^{d-2}} + \frac{d\kappa_d}{2}R^2\Big]. $$
An easy computation (similar to the one in the proof of Lemma \ref{rearr}) gives that
$ \E(H_R^{\{u\}}(x)) = \frac{\kappa_d}{2}|u-x|^2$ if
$|u-x|\le R$.
We have shown above that $U(0) = \lim_{R\to\infty} H_R^{\{0\}}(0)$ converges a.s. Next, by Theorem \ref{thmrearrange}
 it follows that for any $x\in \R^d$ the limit $\lim_{R\to\infty} H_R^{\{0\}}(x)$ exists a.s., uniformly for $x$ in compact sets, since the difference $H_R^{\{0\}}(x)-H_R^{\{0\}}(0)$ can be represented as minus the line integral of $G_R^{\{0\}}(\cdot)$
  (defined in Section \ref{subsectionrearrange}) from $0$ to $x$. By translation, it follows that the limit
$$H^{\{u\}}(x) := \lim_{R\to\infty} H_R^{\{u\}}(x) $$
converges a.s.\ for any fixed $u\in\R^d$ uniformly as $x$ ranges over compact sets. As before,
$H^{\{u\}}(x)$ satisfies a rearrangement identity similar to \eqref{eq:rearrangeidentity}:

\begin{lem}
For any $x,u,v\in\R^d$ we have that a.s.
\begin{equation*}
H^{\{u\}}(x)-H^{\{v\}}(x) = \frac{\kappa_d}{2}\left(|u-x|^2-|v-x|^2\right).
\end{equation*}
\end{lem}

We omit the proof, which is similar to the proof of Lemma \ref{rearr}, and is also superseded by the following stronger lemma.

\begin{lem}
Almost surely, we have that
$$ \max_{u\in B(0,1)} \Big| H_R^{\{u\}}(0)- \frac{\kappa_d}{2}|u|^2 - H_R^{\{0\}}(0)
 \Big| \xrightarrow[R\to\infty]{} 0.$$
\end{lem}

\begin{proof}
For any $m\in \N$ and $\varepsilon>0$, consider the event
$$ J_m^\eps = \Big\{ \max_{m\le R<m+1}
  \max_{u\in B(0,1)} \Big| H_R^{\{u\}}(0)- \frac{\kappa_d}{2}|u|^2 - H_m^{\{0\}}(0)  \Big| > \eps \Big\}.
$$
We shall show that for any $\eps>0$ we have $\sum_{m=1}^\infty \p(J_m^\eps) < \infty$. By Borel-Cantelli, that implies the claim of the lemma.

To that end, fix a large $m\in \N$. Let $E_m$ be a $\frac{1}{m^2}$-net
 of numbers in $[m,m+1]$, and let $N_m$ be a $\frac{1}{m^2}$-net in
 $B(0,1)$, choosing nets such that
 $|E_m\times N_m| = O(m^{2d+2})$ and such that for all $v\in N_m$ and $r\in E_m$ we have
 $ (r-m) \ge m^{-2}/2$ and $|v| \ge m^{-2}/2$. For $v\in N_m$ and $r\in E_m$ denote
$\Delta_{m,v,r} = B(0,m)\triangle B(v,r)$ and
$\nu_{m,v,r} = \text{card}({\cal Z}\cap \Delta_{m,v,r})$ (the number of
 stars in $\Delta_{m,v,r}$). For $v\in N_m, r\in E_m$ define events
\begin{eqnarray*}
K_{m,v,r} &=&
\bigg\{ \Big| \nu_{m,v,r}-\vol(\Delta_{m,v,r}) \Big|
> m^{0.1} \left( \vol(\Delta_{m,v,r})\right)^{1/2}
\bigg\}, \\
L_{m,v,r} &=&
\bigg\{ B\left(v,r+\frac{2}{m^2}\right)\setminus
B\left(v,r-\frac{2}{m^2} \right)\text{ contains }
\\ & & \qquad\qquad\qquad\qquad >20\kappa_d m^{d-3}\text{ stars }
\bigg\},\\
S_{m,v,r}^\eps &=&
\bigg\{ \Big| H_r^{\{v\}}(0) - \frac{\kappa_d}{2}|v|^2 - H_m^{\{0\}}(0) \Big| > \frac{\eps}{2} \bigg\}.
\end{eqnarray*}
Because the number of stars in a region has the Poisson distribution,
by Lemma \ref{poissonlemma} we get that for
some constants $C,c>0$, for all $m$ we have
\begin{equation}\label{eq:desired1}
\p(K_{m,v,r}) \le C e^{-c m^{0.2}}, \qquad \p(L_{m,v,r}) \le C
e^{-c m^{d-3}}
\end{equation}
(note that our choice of the nets $E_m$ and $N_m$ forces $\vol(\Delta_{m,v,r})$
to go to infinity when $m$ grows large).
Next, we derive a bound for $\p(S_{m,v,r}^\eps\cap K_{m,v,r}^c)$. Denote
$$W_{m,v,r} = H_r^{\{v\}}(0) - \frac{\kappa_d}{2}|v|^2 -
H_m^{\{0\}}(0),$$
and observe that $W_{m,v,r}$ is a centered random variable that, conditioned on the event
$\nu_{m,v,r}=k$, can be written as a constant
$ e_{m,v,r} := \frac{1}{d-2}\left(\frac{d\kappa_d}{2}(r^2-m^2)-\frac{\kappa_d}{2}|v|^2\right)$
plus a sum of $k$ i.i.d. random variables with values in $[-(m-1)^{2-d},(m-1)^{2-d}]$.
Therefore we have that
\begin{eqnarray*}
 \E(W_{m,v,r} | \nu_{m,v,r}=k) &=& e_{m,v,r}
+ k\frac{\E(W_{m,v,r}-e_{m,v,r})}{E(\nu_{m,v,r})} \\
&=& e_{m,v,r}\left(1-\frac{k}{\vol(\Delta_{m,v,r})}\right)
\end{eqnarray*}
(since $\E(W_{m,v,r})=0$).
Now, take $k$ such that $$|k-\vol(\Delta_{m,v,r})| \le m^{0.1} \left(\vol(\Delta_{m,v,r})\right)^{1/2}.$$ Noting that for some constant $c_2>0$ we have that
$$\vol(\Delta_{m,v,r}) \ge c_2 m^{d-1}\max(r-m,|v|),$$
it follows that
$$  \E(W_{m,v,r} | \nu_{m,v,r}=k) = m^{0.1} m^{-(d-1)/2} O\left(\frac{m(r-m)+|v|^2}{\max(r-m,|v|)}\right)
=O\left( m^{-0.9} \right)
$$
(since $d\ge 5$). In particular, for such $k$, for any $\eps>0$ we have for sufficiently large $m$ that
$ \E(W_{m,v,r} | \nu_{m,v,r}=k) < \eps/4$, and it follows by Hoeffding's inequality applied to the
representation of $W_{m,v,r}$ conditioned on the event $\{\nu_{m,v,r}=k\}$ described above that for some constants $c_3, c_4>0$ we have that
$$ \p(S_{m,v,r}^\eps\ |\ \nu_{m,v,r}=k) \le e^{-c_3 \eps^2 \frac{m^{2(d-2)}}{\vol(\Delta_{m,v,r})}}
\le e^{-c_4 \eps^2 m^{d-3}}. $$
It follows that
\begin{equation}\label{eq:desired2}
\p(S_{m,v,r}^\eps\cap K_{m,v,r}^c) \le e^{-c_4 \eps^2 m^{d-3}},
\end{equation}
which was our desired estimate.

We now claim that for any fixed $\eps>0$,
for $m$ sufficiently large we have that
$$J_m^\eps \subseteq \bigcup_{v\in N_m, r\in E_m} \Big(K_{m,v,r}\cup
L_{m,v,r} \cup S_{m,v,r}^\eps \Big).
$$
Together with\ \ the above estimates \eqref{eq:desired1} and \eqref{eq:desired2},\ \
 this will prove that $\sum_m \p(J_m^\eps)<\infty$ and therefore the claim of the lemma. To prove this, let $u\in B(0,1)$ and $R\in [m,m+1)$. Let $v\in N_m$ and $r\in E_m$ such that $|v-u|\le m^{-2}$, $|R-r|<m^{-2}$. In particular, we have that the symmetric difference
$ B(u,R) \triangle B(v,r)$ satisfies
$$ B(u,R) \triangle B(v,r) \subseteq B\left(v,r+\frac{2}{m^2}\right)\setminus B\left(v,r-\frac{2}{m^2}\right).
$$
Then
\begin{eqnarray*}
 \Big| H_R^{\{u\}}(0)- \frac{\kappa_d}{2}|u|^2 - H_m^{\{0\}}(0)  \Big| &\le&
 \Big| H_r^{\{v\}}(0)- \frac{\kappa_d}{2}|v|^2 - H_m^{\{0\}}(0)  \Big| \\ & & +
  \Big| H_R^{\{u\}}(0)- \frac{\kappa_d}{2}|u|^2 - H_r^{\{v\}}(0)+\frac{\kappa_d}{2}|v|^2  \Big|
\end{eqnarray*}
Assuming the event $\bigcup_{v\in N_m, r\in E_m}
\Big(K_{m,v,r}\cup L_{m,v,r} \cup S_{m,v,r}^\eps \Big)$ did not occur,
the first term in this bound is $\le \eps/2$, and the second term is at most
$$ \frac{\kappa_d}{2}\Big| |u|^2-|v|^2 \Big| + \frac{20\kappa_d m^{d-3}}{(m/2)^{d-2}}
+ \frac{d\kappa_d}{d-2}\left|R^2-r^2\right|<
\frac{C}{m}. $$
This is also $\le \eps/2$ if $m$ is large enough, which means that $J_m^\eps$ did not occur.
\end{proof}

Combining the above results as in Section \ref{sectionproofprop1}, we have proved:

\begin{prop} \label{potentialdefined}
If $d\ge 5$,
the limit in \eqref{eq:potential} exists a.s.\  simultaneously for all $x \in \R^d\setminus \{z_i\}_i$ and defines a stationary centered process that is a.s.\  differentiable everywhere it is defined and satisfies
\begin{eqnarray}
U(x) &=& \frac{1}{d-2}
\lim_{T\to\infty} \Big[\sum_{i\,:\,|z_i|<T} \frac{-1}{|z_i-x|^{d-2}} + \frac{d\kappa_d}{2}T^2\Big]
- \frac{\kappa_d}{2}|x|^2, \label{eq:ualtern} \\
\nabla U(x) &=& -F(x). \label{eq:gradu}
\end{eqnarray}
\end{prop}

We will occasionally use truncated versions of the gravitational
potential, the force and its first differential. For a bounded set
$A\subset \R^d$, define $U(x | A)$,
the {\bf partial potential from stars in $A$}, by
$$U(x | A)= \frac{1}{d-2}\sum_{i\,:\,z_i\in A} \frac{-1}{|z_i-x|^{d-2}} +
\frac{1}{d-2}\int_A |z-x|^{-d+2} d\vol(z). $$
Similarly, define $F(x | A)$, the {\bf partial force from stars in
$A$} by
$$
F(x | A) = \sum_{i\,:\,z_i\in A} \frac{z_i-x}{|z_i-x|^{d}}
  - \int_A \frac{z-x}{|z-x|^d} d\vol(z). $$
For a set $A \subset \R^d$ whose complement is bounded, define
\begin{eqnarray*}
U(x | A)&=& U(x) - U(x | A^c), \\
F(x | A) &=& F(x) - F(x | A^c).
\end{eqnarray*}
Note that $U(x | A)$ and $F(x | A)$ are centered to have mean 0.

While these definitions are rather general, throughout the paper we
only use sets $A$ which are annuli
of the form $A=B(y,p)\setminus B(y,q)$, where $0\le q<p\le \infty$
(this includes the degenerate cases of a ball, the complement of a
ball, and the entire space). Furthermore, in all the cases we will
consider, we will have that either [$q>0$ and $|x-y|\le q$] or [$q=0$ and
$|x-y|\le p$]. In those cases, from
the computations in the proof above it is easy to verify that we
have the following explicit expressions for $U(x | A)$ and $F(x | A)$:
First, if $q>0$ and $|x-y| \le q$, then
\begin{eqnarray} \nonumber
U(x | B(y,p)\setminus B(y,q)) &=& \frac{1}{d-2}\sum_{q<|z_i-y|\le p,
\ |z_i| \uparrow}
\frac{-1}{|z_i-x|^{d-2}} \\ & &+ \frac{d \kappa_d}{2(d-2)}(p^2-q^2),
\label{eq:numberedeq1} \\
F(x | B(y,p)\setminus B(y,q)) &=& \sum_{q<|z_i-y|\le p,\ |z_i|\uparrow}
\frac{z_i-x}{|z_i-x|^d}. \nonumber
\end{eqnarray}
Second, if $q=0$ and $|x-y|\le p$ then
\begin{eqnarray*}
U(x | B(y,p)) &=& \frac{1}{d-2}\sum_{|z_i-y|\le p,\ |z_i|\uparrow}
\frac{-1}{|z_i-x|^{d-2}} + \frac{d\kappa_d}{2(d-2)}p^2, \\
F(x | B(y,p)) &=& \sum_{|z_i-y|\le p,\ |z_i|\uparrow}
\frac{z_i-x}{|z_i-x|^d} + \kappa_d(x-y).
\end{eqnarray*}
We will also use the function $D_1 F(x | A)$, the first differential
of $F(x | A)$. By Lemma \ref{simconv1} we have the following explicit expressions
for
$D_1 F(x | A)$ in the cases described above: If $q>0$ and $|x-y|\le q$ then
$$ D_1 F(x | B(y,p)\setminus B(y,q)) = \sum_{q<|z_i-y|\le p,\
|z_i|\uparrow} D_1\left[ \frac{z_i-x}{|z_i-x|^d} \right], $$
and if $q=0$ and $|x-y|\le p$ then
\begin{equation}\label{eq:numberedeq2}
 D_1 F(x | B(y,p)) = \sum_{|z_i-y|\le p,\
|z_i|\uparrow} D_1\left[ \frac{z_i-x}{|z_i-x|^d} \right] + \kappa_d
I_{d\times d},
\end{equation}
where $I_{d\times d}$ is the $d\times d$ identity operator.

In all the above sums, if the region of summation is infinite then the
terms are summed in order of increasing distance from $0$.

\section{Large deviations estimates \label{sectionlarge}}

In this section we derive detailed large deviations estimates for the force $F(x)$, its derivative $D_1 F$, and (in dimensions 5 and higher) the gravitational potential function.

\subsection{Large deviations for the value at a point \label{subsectionlargepoint}}

Consistently with the previously defined notation, let $D_k F\big(x\ \big|\ A\big)$ denote the $k$-th differential tensor of the function $x\to F\big(x\ \big|\ A\big)$.

\begin{thm} \label{largedeviations1} There exist constants $C_1, c_2, c_3>0$ such that
for all $p > q > 0$ and $t>0$ we have
\begin{eqnarray}
\p\Big(\Big|U\big(0\ \big|\ B(0,p) \setminus B(0,q)\big)\Big|\ge t \Big) &\le&
C_1 e^{-c_2 q^{d-2} t \log\left(\frac{c_3 t}{q^2}\right)}, \label{eq:firstone} \\
\p\Big(\Big|F\big(0\ \big|\ B(0,p) \setminus B(0,q)\big)\Big|\ge t \Big) &\le&
C_1 e^{-c_2 q^{d-1} t \log\left(\frac{c_3 t}{q}\right)}, \label{eq:secondone} \\
\p\Big(\Big|D_1 F\big(0\ \big|\ B(0,p) \setminus B(0,q)\big)\Big|\ge t \Big) &\le&
C_1 e^{-c_2 q^d t \log\left(c_3 t \right)}. \label{eq:thirdone}
\end{eqnarray}
Equation \eqref{eq:firstone} holds in dimensions $d\ge 5$, and equations \eqref{eq:secondone} and \eqref{eq:thirdone} hold for all dimensions $d\ge 3$.
\end{thm}

\begin{proof} Assume $d\ge 5$.
Denote $B_{p,q} = B(0,p)\setminus B(0,q)$. Let $W_{p,q}, N_{p,q}, U_{p,q}$ be as in Section \ref{sectionpotential}, so
$U\big(0\ \big|\ B(0,p) \setminus B(0,q)\big)  = \frac{1}{d-2}(U_{p,q}-\E[U_{p,q}])$.
Let $V=|W_{p,q}|^{-(d-2)}$. Then for any $u>0$ we have
\begin{eqnarray*}
\p(|V|>u) &=& \p\left(\frac{1}{|W_{p,q}|^{d-2}}>u\right) = \p\left( |W_{p,q}|< \frac{1}{u^{1/(d-2)}} \right)
\\ &\le& \frac{\kappa_d}{\vol(B_{p,q}) u^{d/(d-2)}}.
\end{eqnarray*}
Therefore, noting that $0\le V \le 1/q^{d-2}$, we have for any integer $k\ge 2$ that
\begin{eqnarray}
\E(|V|^k)  &=& \int_0^{1/q^{d-2}} k u^{k-1}\p(|V|>u)du \le \frac{\kappa_d}{\vol(B_{p,q})}
\int_0^{1/q^{d-2}} k u^{k-2-\frac{2}{d-2}} du \nonumber \\ &=&
\frac{\kappa_d}{\vol(B_{p,q})}\cdot
\frac{k}{k-1-\frac{2}{d-2}} \left(\frac{1}{q^{d-2}}\right)^{k-1-\frac{2}{d-2}} \nonumber \\ &\le&
\frac{6\kappa_d}{\vol(B_{p,q}) q^{(d-2)k-d}}. \label{eq:moment}
\end{eqnarray}
For any $\theta\ge0$ we have
\begin{eqnarray*}
\E(e^{\theta V}) &=& 1 + \theta\E(V)+\sum_{k=2}^\infty \frac{\theta^k}{k!}\E V^k \le
1+ \theta\E(V)+ \frac{6\kappa_d q^d}{\vol(B_{p,q})} \sum_{k=2}^\infty \frac{(\theta/q^{d-2})^k}{k!}
\\&\le& 1+ \theta\E(V) + \frac{6\kappa_d q^d}{\vol(B_{p,q})} e^{\theta/q^{d-2}}.
\end{eqnarray*}
Conditionally on $N_{p,q}$, the stars in $B_{p,q}$ are a vector of $N_{p,q}$ i.i.d points distributed uniformly in $B_{p,q}$. Hence by the last estimate we get that
$$ \E\Big[e^{\theta U_{p,q}}\ \big|\ N_{p,q}\Big]
\le \left(1+ \theta\E(V) + \frac{6\kappa_d q^d}{\vol(B_{p,q})} e^{\theta/q^{d-2}}\right)^{N_{p,q}}. $$
Now, it is a simple exercise that if $X\sim \text{Poi}(\lambda)$, then for any $\alpha$, we have $\E((1+\alpha)^X)=e^{\lambda \alpha}$. Since $N_{p,q} \sim \text{Poi}(\vol(B_{p,q}))$,
using this in the above inequality we get
\begin{equation}\label{eq:weget}
\E\left(e^{\theta U_{p,q}}\right) \le \exp\Big(\theta\vol(B_{p,q})\E(V)+ 6\kappa_d q^d e^{\theta/q^{d-2}}\Big).
\end{equation}
Since also $\E[U_{p,q}]=\vol(B_{p,q})\E(V)$, we get
$$
\E\left(e^{\theta (U_{p,q}-\E(U_{p,q}))}\right) \le \exp\Big(6\kappa_d q^d e^{\theta/q^{d-2}}\Big).
$$
Therefore for any $t\ge 0$ we have
$$ \p\Big(U_{p,q}-\E[U_{p,q}] \ge t\Big) \le \exp\left( - \theta t + 6\kappa_d q^d e^{\theta/q^{d-2}} \right). $$
Set $\theta=q^{d-2} \log(t/(6\kappa_d q^2))$ to get
\begin{eqnarray*}
\p\Big(U_{p,q}-\E[U_{p,q}]\ge t\Big) &\le& \exp\Big(- q^{d-2} t \log(t/(6\kappa_d q^2))+q^{d-2} t\Big)
\\ &=& \exp\bigg(-q^{d-2} t \log \left(\frac{t}{6e \kappa_d q^2}\right)  \bigg).
\end{eqnarray*}
In the same way, one gets a similar bound for the negative tail $\p(U_{p,q}-\E(U_{p,q}) < -t)$, by noting that \eqref{eq:weget} also holds for negative values
of $\theta$ if on the right-hand side $e^{\theta/q^{d-2}}$ is
replaced by $e^{|\theta|/q^{d-2}}$.
Combining the negative and positive tail bounds gives
\eqref{eq:firstone}. The estimates \eqref{eq:secondone} and \eqref{eq:thirdone} follow (with the weaker assumption $d\ge 3$) by estimating in exactly the same way the moments and exponential moments of $|W_{p,q}|^{-(d-1)}$ and $|W_{p,q}|^{-d}$, respectively, in place of $|W_{p,q}|^{-(d-2)}$. Note that the random variables
$U\big(0\ \big|\ B(0,p) \setminus B(0,q)\big),
F\big(0\ \big|\ B(0,p) \setminus B(0,q)\big),
D_1 F\big(0\ \big|\ B(0,p) \setminus B(0,q)\big)$ are all centered. We omit the full proofs.
\end{proof}

\subsection{Uniform bounds in a ball \label{subsectionlargeball}}

\begin{thm} \label{largedeviations2} There exist constants $C_1, c_2, c_3>0$ such that
for all $p > q > 0$ and $t>0$ we have
\begin{eqnarray}
\p\Big(\max_{x\in B(0,1\wedge \frac{q}{2})} \Big|U\big(x\ \big|\ B(0,p) \setminus B(0,q)\big)\Big|\ge t \Big) &\le&
C_1 e^{-c_2 q^{d-2} t \log\left(\frac{c_3 t}{q^2}\right)}, \label{eq:ballfirstone} \\
\p\Big(\max_{x\in B(0,1\wedge \frac{q}{2})}\Big|F\big(x\ \big|\ B(0,p) \setminus B(0,q)\big)\Big|\ge t \Big) &\le&
C_1 e^{-c_2 q^{d-1} t \log\left(\frac{c_3 t}{q}\right)}, \label{eq:ballsecondone} \\
\p\Big(\max_{x\in B(0,1\wedge \frac{q}{2})} \Big|D_1 F\big(x\ \big|\ B(0,p) \setminus B(0,q)\big)\Big|\ge t \Big) &\le&
C_1 e^{-c_2 q^d t \log\left(c_3 t \right)}, \label{eq:ballthirdone}
\end{eqnarray}
where equation \eqref{eq:ballfirstone} holds in dimensions $d\ge 5$, and equations \eqref{eq:ballsecondone} and \eqref{eq:ballthirdone} hold for all dimensions $d\ge 3$.

\end{thm}

\begin{proof} Start with \eqref{eq:ballthirdone}. Set $r=t/q$. We have

$$ \p\Big(\max_{x\in B(0,1\wedge \frac{q}{2})} \Big|D_1 F\big(x\ \big|\ B(0,p) \setminus B(0,q)\big)\Big|\ge t \Big) \qquad\qquad\qquad\qquad$$

\vspace{-20.0pt}
\begin{eqnarray} \nonumber
 &\le&
\p\Big(\Big|D_1 F\big(0\ \big|\ B(0,p) \setminus B(0,q)\big)\Big|\ge \frac{t}{2} \Big)
\\ & & +
\p\Big(\max_{x\in B(0,1\wedge \frac{q}{2})} \Big|D_2 F\big(x\ \big|\ B(0,p) \setminus B(0,q)\big)\Big|\ge
r \Big). \label{eq:twoline}
\end{eqnarray}
The first of these two terms is bounded by $C_1 e^{-c_2 q^d t\log(c_3 t)}$ by \eqref{eq:thirdone}. For the second term, observe that by \eqref{eq:apriori2} we have that
$$
\p\Big(\max_{x\in B(0,1\wedge \frac{q}{2})}
\Big|D_2 F\big(x\ \big|\ B(0,p) \setminus B(0,q)\big)\Big|\ge r \Big)
\qquad\qquad\qquad\qquad\qquad
$$

\vspace{-17.0pt}
\begin{eqnarray*} &\le&
\sum_{m=0}^\infty \p\left(\max_{x\in
B\left(0,1\wedge\frac{q}{2}\right)} \Big|D_2 F\big(x\ \big|\ B(0,p)\cap(B(0,2^{m+1}
q)\setminus B(0,2^m q))\big)\Big|\ge \frac{r}{2^{m+1}}\right)
\\ &\le&
\sum_{m=0}^\infty \p\left(\frac{C \nu_m}{2^{m(d+1)}q^{d+1}} \ge
\frac{r}{2^{m+1}}\right),
\end{eqnarray*}
where $\nu_m$ is the number of stars in $B(0,p)\cap
(B(0,2^{m+1}q)\setminus B(0,2^m q))$, which is a Poisson random variable
with mean $\le C2^{dm}q^d$.
Using Lemma \ref{poissonlemma} it follows that the above sum is less than
$$ \sum_{m=0}^\infty C \exp\left(-c r 2^{md}q^{d+1}\log(c q
r)\right). $$
Now, if in the above inequality $crq^{d+1}\log(cqr)>2$ then the whole
sum is less than a constant times its first term, so
$$ \p\left(\max_{x\in B\left(0,1\wedge \frac{q}{2}\right)}
\Big|D_2 F\big(x\ \big|\ B(0,p)\setminus B(0,q)\big)\Big| \ge r \right)
\le C \exp\left(-crq^{d+1}\log(c q r)\right). $$
On the other hand, if $crq^{d+1}\log(cqr)\le 2$ then the above
inequality holds trivially if we take $C$ slightly larger since then
the right-hand side is larger than $1$. Hence this inequality holds
for all values of $r$ and $q$. Plugging this into equation
\eqref{eq:twoline} together with the fact mentioned after \eqref{eq:twoline}
gives \eqref{eq:ballthirdone}.

Next, to prove \eqref{eq:ballsecondone}, write
$$ \p\Big(\max_{x\in B(0,1\wedge \frac{q}{2})} \Big|F\big(x\ \big|\ B(0,p) \setminus B(0,q)\big)\Big|\ge t \Big) \qquad\qquad\qquad\qquad$$

\vspace{-20.0pt}
\begin{eqnarray*}
 &\le&
\p\Big(\Big|F\big(0\ \big|\ B(0,p) \setminus B(0,q)\big)\Big|\ge \frac{t}{2} \Big)
\\ & & +
\p\Big(\max_{x\in B(0,1\wedge \frac{q}{2})} \Big|D_1 F\big(x\ \big|\
B(0,p) \setminus B(0,q)\big)\Big|\ge r \Big),
\end{eqnarray*}
where again $r=t/q$.
Both of the terms are bounded by $C_1 e^{-c_2 q^{d-1} t\log(c_3 t/q)}$ by \eqref{eq:secondone} and \eqref{eq:ballthirdone}.

Finally, to prove \eqref{eq:ballfirstone}, write similarly
$$ \p\Big(\max_{x\in B(0,1\wedge \frac{q}{2})} \Big|U\big(x\ \big|\ B(0,p) \setminus B(0,q)\big)\Big|\ge t \Big) \qquad\qquad\qquad\qquad\qquad\qquad\qquad$$

\vspace{-20.0pt}
\begin{eqnarray*}
 &\le&
\p\Big(\Big|U\big(0\ \big|\ B(0,p) \setminus B(0,q)\big)\Big|\ge \frac{t}{2} \Big)
\\ & & +
\p\Big(\max_{x\in B(0,1\wedge \frac{q}{2})} \Big|F\big(x\ \big|\ B(0,p) \setminus B(0,q)\big)\Big|\ge
r \Big) \le C_1 e^{-c_2 q^{d-2} t\log(c_3 t/q^2)}
\end{eqnarray*}
by \eqref{eq:firstone} and \eqref{eq:ballsecondone}.
\end{proof}

By letting $p\to\infty$ we get the following limiting case of Theorem \ref{largedeviations2}.

\begin{cor} \label{limlargedeviations2} There exist constants $C_1, c_2, c_3>0$ such that
for all $q > 0$ and $t>0$ we have
\begin{eqnarray}
\p\Big(\max_{x\in B(0,1\wedge \frac{q}{2})} \Big|U\big(x\ \big|\ \R^d \setminus B(0,q)\big)\Big|\ge t \Big) &\le&
C_1 e^{-c_2 q^{d-2} t \log\left(\frac{c_3 t}{q^2}\right)}, \label{eq:ballfirstonelim} \\
\p\Big(\max_{x\in B(0,1\wedge \frac{q}{2})}\Big|F\big(x\ \big|\ \R^d \setminus B(0,q)\big)\Big|\ge t \Big) &\le&
C_1 e^{-c_2 q^{d-1} t \log\left(\frac{c_3 t}{q}\right)}, \label{eq:ballsecondonelim} \\
\p\Big(\max_{x\in B(0,1\wedge \frac{q}{2})} \Big|D_1 F\big(x\ \big|\ \R^d \setminus B(0,q)\big)\Big|\ge t \Big) &\le&
C_1 e^{-c_2 q^d t \log\left(c_3 t \right)}, \label{eq:ballthirdonelim}
\end{eqnarray}
where equation \eqref{eq:ballfirstonelim} holds in dimensions $d\ge 5$, and equations \eqref{eq:ballsecondonelim} and \eqref{eq:ballthirdonelim} hold for all dimensions $d\ge 3$.
\end{cor}


\subsection{Uniform bounds in a ball with a moving domain \label{subsectionlargemov}}

\begin{thm} \label{largedeviations3} There exist constants $C_1, c_2, c_3, C_4>0$ such that
for all $p > q > 0$ and $t>0$ we have that if $t>C_4 p^2$ or $p=\infty$ then
\begin{eqnarray}
\p\Big(\max_{x\in B(0,1\wedge \frac{q}{2})} \Big|U\big(x\ \big|\ B(x,p) \setminus B(x,q)\big)\Big|\ge t \Big) &\le&
C_1 e^{-c_2 q^{d-2} t \log\left(\frac{c_3 t}{q^2}\right)}, \label{eq:ballfirstonemov} \\
\p\Big(\max_{x\in B(0,1\wedge \frac{q}{2})}\Big|F\big(x\ \big|\ B(x,p) \setminus B(x,q)\big)\Big|\ge t \Big) &\le&
C_1 e^{-c_2 q^{d-1} t \log\left(\frac{c_3 t}{q}\right)}, \label{eq:ballsecondonemov} \\
\p\Big(\max_{x\in B(0,1\wedge \frac{q}{2})} \Big|D_1 F\big(x\ \big|\ B(x,p) \setminus B(x,q)\big)\Big|\ge t \Big) &\le&
C_1 e^{-c_2 q^d t \log\left(c_3 t \right)}, \label{eq:ballthirdonemov}
\end{eqnarray}
where equation \eqref{eq:ballfirstonemov} holds in dimensions $d\ge 5$, and equations \eqref{eq:ballsecondonemov} and \eqref{eq:ballthirdonemov} hold for all dimensions $d\ge 3$.
\end{thm}

\begin{proof}
Denote $B=B\left(0,1\wedge \frac{q}{2}\right)$. First, we prove
\eqref{eq:ballfirstonemov} in the limiting case when $p=\infty$. Fix
$x\in B$, then
$$
\Big|U\big(x\ \big|\ \R^d \setminus B(x,q)\big) -
U\big(x\ \big|\ \R^d \setminus B(0,q)\big)\Big|
\qquad\qquad\qquad\qquad $$

\vspace{-18.0pt}
\begin{eqnarray*} \qquad\qquad &=&
\Big|U\big(x\ \big|\ B(0,q)\big) -
U\big(x\ \big|\ B(x,q)\big)\Big| \\ &=&
-\kappa_d |x|^2/2 - \frac{1}{d-2}\sum_{z_i\in E_1}
\frac{1}{|z_i-x|^{d-2}} + \frac{1}{d-2}\sum_{z_i\in E_2}
\frac{1}{|z_i-x|^{d-2}},
\end{eqnarray*}
where $E_1 = B(0,q)\setminus B(x,q)$ and $E_2 = B(x,q)\setminus
B(0,q)$. Now, denoting by $\nu_q$ the number of stars in
$B\left(0,q+1\wedge
\frac{q}{2}\right)-B\left(0,q-1\wedge\frac{q}{2}\right)$, it follows
that
$$
\Big|U\big(x\ \big|\ \R^d \setminus B(x,q)\big)\Big| \le
\Big|U\big(x\ \big|\ \R^d \setminus B(0,q)\big)\Big| + C_5 q^2 +
\frac{\nu_q}{(d-2)(q/2)^{d-2}}. $$
Since $\nu_q$ is a Poisson random
variable with mean $\le C q^d$, by Lemma \ref{poissonlemma} we obtain that for
$t>3C_5 q^2$ we have
\begin{eqnarray}
\p\left( \max_{x\in B}
\Big|U\big(x\ \big|\ \R^d \setminus B(x,q)\big)\Big| \ge t \right)
&\le&
\p\left( \max_{x\in B}
\Big|U\big(x\ \big|\ \R^d \setminus B(0,q)\big)\Big| \ge t/3 \right)
\nonumber \\ \nonumber & & +
\p\left( \frac{\nu_q}{(d-2)(q/2)^{d-2}}\ge t/3 \right) \\ &\le& C
\exp\left(-ctq^{d-2} \log\left(\frac{ct}{q^2}\right)\right).
\label{eq:haveproved}
\end{eqnarray}
This also holds trivially for $t\le 3C_5 q^2$ (since the RHS is larger
than $1$) provided $c$ is chosen small enough, hence it gives
\eqref{eq:ballfirstonemov} in the case $p=\infty$.
%
%
%
%

Finally, to prove \eqref{eq:ballfirstonemov} in the general case, note, using \eqref{eq:haveproved} twice
and using the assumption $t>C_4 p^2$, that
$$
\p\Big(\max_{x\in B(0,1\wedge \frac{q}{2})} \Big|U\big(x\ \big|\ B(x,p) \setminus B(x,q)\big)\Big|\ge t \Big)
\qquad\qquad\qquad\qquad\qquad
$$

\vspace{-20.0pt}
\begin{eqnarray*} &\le &
\p\Big(\max_{x\in B(0,1\wedge \frac{q}{2})} \Big|U\big(x\ \big|\ \R^d \setminus B(x,q)\big)\Big|\ge
\frac{t}{2} \Big)
\\ & & +
\p\Big(\max_{x\in B(0,1\wedge \frac{q}{2})} \Big|U\big(x\ \big|\ \R^d \setminus B(x,p)\big)\Big|\ge
\frac{t}{2} \Big)
\\ & \le &
C_1 e^{-c_2 q^{d-2} t \log\left(\frac{c_3 t}{q^2}\right)} +
C_1 e^{-c_2 p^{d-2} t \log\left(\frac{c_3 t}{p^2}\right)}.
\\ & \le &
2C_1 e^{-c_2 q^{d-2} t \log\left(\frac{c_3 t}{q^2}\right)}.
\end{eqnarray*}
The proofs of \eqref{eq:ballsecondonemov} and \eqref{eq:ballthirdonemov} are similar and are omitted.
\end{proof}

\section{Proof of Theorem \ref{thm3} in dimensions 5 and higher
\label{sectionproofthm3}}

In this section, we assume that $d\ge 5$. Our goal is to bound the probability of the event $E_R$ of having a gravitational flow curve connect $\partial Q(0,R)$ with $\partial Q(0,2R)$, as $R\to\infty$.
The case of dimensions 3 and 4 is slightly more delicate. In Section \ref{sectiondim34} we explain what modifications to the proof are required to complete the proof of Theorem \ref{thm3} in that case.

\subsection{Reduction to a problem on a discrete set of points \label{subsectiondiscrete}}

Fix the following parameters:
\begin{eqnarray*}
B &=& R^{8/9}, \\
\Delta &=& \textrm{a large constant (depending on $d$) whose value will be specified later},\\
r &=& \Delta\cdot (\log R)^{2/d}, \\
\rho &=& R^{-1/10}, \\
s &=& R^{-\frac{1}{10(d^2+1)}}, \\
\varepsilon &=& \frac{\rho}{s^d}\log R.
\end{eqnarray*}
We emphasize that $R$ is the only true parameter here, and the values of all the other quantities are specified as functions of $R$.

To control the event $E_R$, we discretize space. Introduce a grid of points in the region $Q(0,2R)\setminus Q(0,R)$, defined by
$$ S = r\Z^d \cap (Q(0,2R)\setminus Q(0,R)). $$
We think of $S$ as an induced subgraph of $r\Z^d$ with the usual lattice structure. Thus, two points
$w,w'\in S$ are called {\bf adjacent} if $|w-w'|=r$. A set $W\subset S$ is called {\bf connected} if
the induced subgraph of $W$ in $S$ is connected. Call a set $W\subset S$ {\bf connectible} if
$W$ is contained in a set $W' \subset S$ which is connected and $|W'|\le 10^d |W|$.
To each point $w\in S$ associate an {\bf inner box} $Q_{\textrm{in}}(w) = Q(w,r)$ and an {\bf outer box} $Q_{\textrm{out}}(w) = Q(w,2r)$.

\begin{lem}\label{fewconnectible}
There exists a constant $C_{15}>0$ such that for any $L\ge 1$
the number of connectible sets $W\subset S$ of cardinality $L$ is at most $R^d C_{15}^L$.
\end{lem}

\begin{proof}
This is an immediate consequence of \cite[Eq. (4.24), p. 81]{grimmett}
\end{proof}

\begin{defi} \label{defibad}
Say that $w\in S$ is {\bf bad} if there exists a gravitational flow curve $\gamma$ connecting $\partial Q_{\textrm{in}}(w)$ with $\partial Q_{\textrm{out}}(w)$ such that at least one of the following conditions hold: \\
(1) $U\big(x\ |\ B(x,3R^{1/d})\big) < -\frac{B}{2} $ for all $x\in \gamma$, or \\
(2) $ \int_\gamma |F(x)|\cdot |dx| < \frac{10^d B r}{R}. $
\end{defi}

We wish to show that the ``bad'' event $E_R$, whose probability we are trying to bound, implies the occurrence of many bad grid points. This will be true except on some atypical events which will happen with probability small enough as to be of no consequence. Define
\begin{eqnarray*}
\Omega_1 &=& \Big\{ \max_{x\in Q(0,2R)} U(x)>B \Big\}, \\
\Omega_2 &=&
  \Big\{ \max_{x\in Q(0,2R)} \Big|U\big(x\ \big|\  \R^d\setminus B(x,3 R^{1/d})\big)\Big| \ge \frac{B}{2} \Big\}, \\
\Omega_3 &=&
  \Big\{ \max_{x\in Q(0,3R)}
    \Big|D_1 F\big(x\ \big|\  \R^d\setminus B(x,3 R^{1/d})\big)\Big| \ge \frac{\varepsilon}{8\sqrt{d}\rho} \Big\}.
\end{eqnarray*}

\begin{lem} \label{exceptionalomega}  For some constants $C,c>0$ we have for all $R$ sufficiently
large that
\begin{eqnarray}
\p(\Omega_1) &\le& C e^{-c R^{4/3}}, \label{eq:omega1} \\
\p(\Omega_2) &\le& C e^{-c R^{11/9}} \label{eq:omega2}, \\
 \p(\Omega_3) &\le&C e^{-c R^{1+1/100 d^2}} \label{eq:omega25}.
\end{eqnarray}
\end{lem}

\begin{proof} First, we prove \eqref{eq:omega2}.
Cover $Q(0,2R)$ with $O(R^d)$ balls $\{B_j\}_{j=1}^{C R^d}$ of radius $1$.
For each ball $B_j$ we have by Theorem \ref{largedeviations3} that
$$ \p\bigg(\max_{x\in B_j} \Big|U\big(x\ \big|\  \R^d\setminus B(x,3R^{1/d})\big)\Big|\ge\frac{B}{2} \bigg)
\le C e^{-c R^{(d-2)/d} B \log(c B/R^{2/d})}. $$
Therefore by a union bound we get that for some new constant $C'>0$, 
\begin{equation*}
\p(\Omega_2) \le
C' R^d e^{-c R^{(d-2)/d} B \log(c B/R^{2/d})}.
\end{equation*}
Now substitute the values of the parameters to get \eqref{eq:omega2}.

Next, \eqref{eq:omega25} follows from Theorem \ref{largedeviations3} in the same way as
\eqref{eq:omega2}.

Finally, to prove \eqref{eq:omega1},
let $a>0$ be some small positive number such that $a<((d-2)/d\kappa_d)^{1/2}$ (another condition will be imposed on it shortly). Note that
$$ U(x) = U\big(x\ \big|\ B(x,a \sqrt{B})\big) + U\big(x\ \big|\ \R^d\setminus B(x,a \sqrt{B})\big),
$$
and that $U\big(x\ \big|\ B(x,a \sqrt{B})\big) \le \frac{d\kappa_d}{2(d-2)} a^2 B < \frac{B}{2}$ (see eq. \eqref{eq:numberedeq1}), so that
on $\Omega_1$ we have that
$$ \max_{x\in Q(0,2R)} \Big|U\big(x\ \big|\ \R^d\setminus B(x,a \sqrt{B})\big)\Big| > \frac{B}{2}. $$
By a similar argument to that used in the proof of \eqref{eq:omega2} above,
the probability of this is bounded by $C R^d e^{-ca^{d-2} B^{d/2}\log(c/a^2)}$. If $a$ was chosen sufficiently small this gives the bound \eqref{eq:omega1} upon substituting the values of the parameters.
\end{proof}

\begin{lem} \label{manybad}
On the event $\Omega_1^c\cap \Omega_2^c$, if there exists a gravitational flow curve $\Gamma$ connecting $\partial Q(0,R)$ and $\partial Q(0,2R)$ (that is, if $E_R$ occurred), and if $R$ is large enough, then
there exists a connectible family $W \subseteq S$ of bad points, with $|W| \ge R/10^d r$.
\end{lem}

\begin{figure}
\centering \resizebox{!}{7cm}{\includegraphics{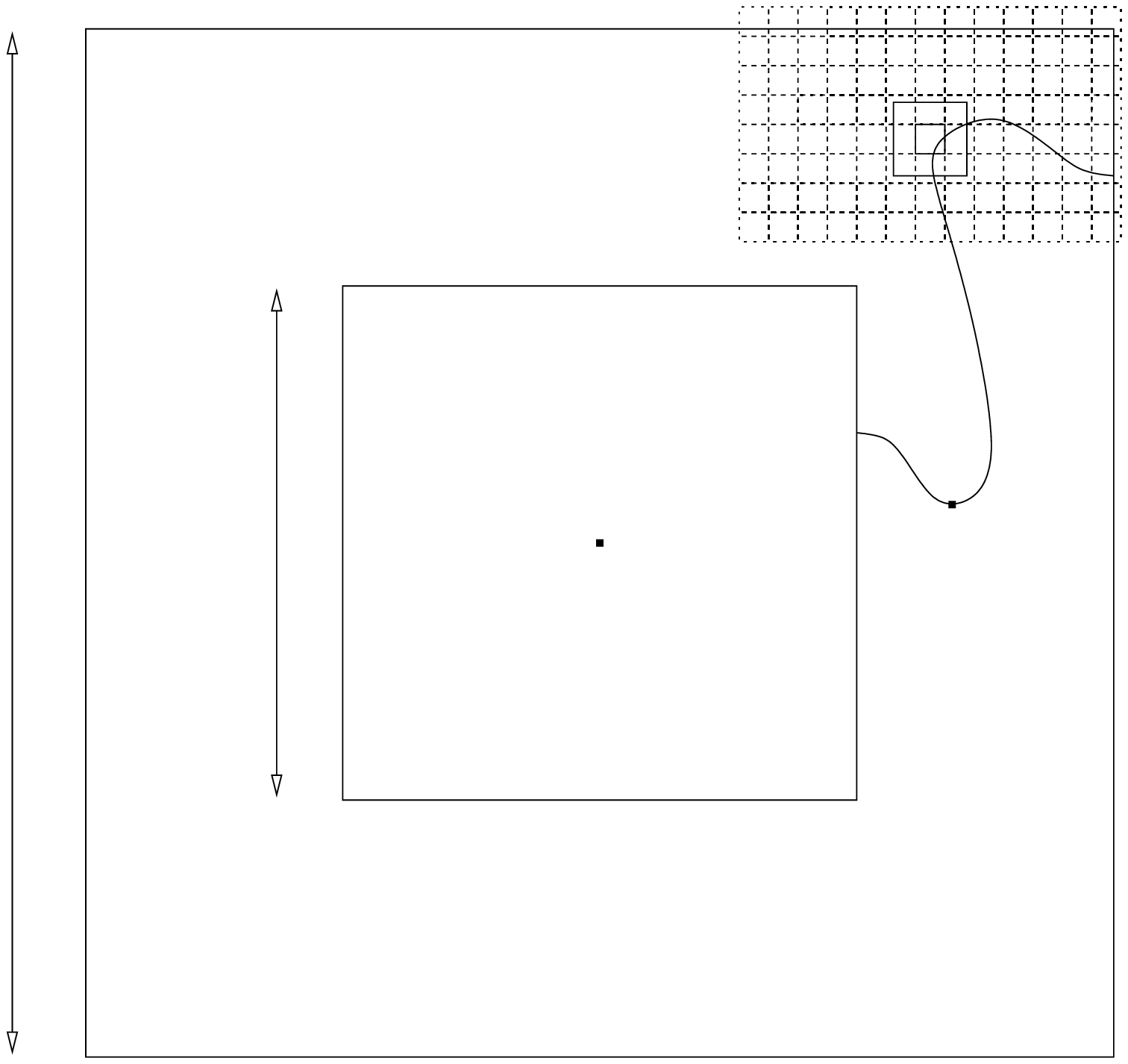}}

\vspace{-110.0pt} \hspace{-236.0pt} {\small $4R$}

\vspace{-15.0pt} \hspace{-134.0pt} {\small $2R$}

\vspace{-15.0pt} \hspace{-8.0pt} {\small $0$}

\vspace{-12.0pt} \hspace{145.0pt} {\small $\Gamma(T_1)$}

\vspace{-100.0pt} \hspace{280.0pt} {\small $\bigg\}$ The grid $S$}

\vspace{30.0pt} \hspace{170.0pt} {\small $\Gamma$}

\vspace{130.0pt}

\caption{Schematic illustration of the proof of Lemma \ref{manybad}. \label{gridfig}}
\end{figure}

\begin{proof} Let $\Gamma:[0,T]\to \R^d$ be a flow curve that connects $\partial Q(0,R)$ and $\partial Q(0,2R)$, and assume that $\Omega_1\cup \Omega_2$ did not occur. In particular,
$U(\Gamma(0))\le B$. Observe that the potential $U$ decreases along the curve $\Gamma$, since $F(x) = - \nabla U(x)$, and therefore
$$ \frac{d}{dt} U(\Gamma(t)) = \langle \dot{\Gamma}(t), \nabla
U(\Gamma(t)) \rangle =
\langle F(\Gamma(t)), -F(\Gamma(t)) \rangle = -|F(\Gamma(t))|^2. $$

Let $W'$ be the set of points $w\in S$ such that $\Gamma$ intersects both $\partial Q_{\textrm{in}}(w)$ and $\partial Q_{\textrm{out}}(w)$.
Since $\Gamma$ connects $\partial Q(0,R)$ and $\partial Q(0,2R)$, clearly we have that $|W'|\ge R/r-2$ (the $-2$ is to account for boundary effects).

Let $T_1 \in [0,T]$ be the least time for which $U(\Gamma(T_1))\le  -B$, or let
$T_1=T$ if no such time exists (see Figure \ref{gridfig}). Certainly,
all the points $w\in S$ for which $\Gamma_{\big| [T_1,T]}$ intersects both $\partial Q_{\textrm{in}}(w)$ and
$\partial Q_{\textrm{out}}(w)$ are bad (since, because $\Omega_2^c$ occurred, they satisfy condition (1) in the definition).
If there are $R/10^d r$ such points, we are done, since the set of such $w$ is connected and \emph{a fortiori} connectible. If this is not so, denote by $W''$ the set of those $w \in W'$ for which
$\Gamma_{\big| [0,T_1]}$ intersects both $\partial Q_{\textrm{in}}(w)$ and $\partial Q_{\textrm{out}}(w)$. The family $W''$ is
a connected set, and we have $|W''| \ge |W'| - R/10^d r - 4^d > R/2r$ (for $R$ large; the $4^d$ is
again to account for boundary effects near $\Gamma(T_1)$).
For each $w\in W''$ let $\Gamma_w$ denote some segment of $\Gamma_{\big|[0,T_1]}$ that connects $\partial Q_{\textrm{in}}(w)$ with $\partial Q_{\textrm{out}}(w)$ (possibly in the opposite direction) and that is contained in the interior of $Q_{\text{out}}(w)\setminus Q_{\text{in}}(w)$ except for its endpoints. Note that the segments
$(\Gamma_w)_{w\in W''}$ are not necessarily disjoint. Replace $W''$ by a subset $W'''\subset W''$ such that the interiors of $(Q_{\textrm{out}}(w))_{w\in W'''}$ are disjoint (and therefore also $(\Gamma_w)_{w\in W'''}$ are disjoint except possibly for their endpoints) and $|W'''|\ge |W''|/5^d$. This can be done using a greedy method, since each point $w\in W''$ added to $W'''$ eliminates at most $5^d$ others.

Let $k$ denote the number of $w\in W'''$ which are not bad. Then
\begin{eqnarray*}
2B &\ge& U(\Gamma(0)) - U(\Gamma(T_1)) \\ &=&
 \int_{\Gamma_{\big| [0,T_1]}} |F(x)|\cdot |dx| \ge
\int_{\bigcup_{w\in W'''} \Gamma_w} |F(x)|\cdot |dx| \ge k \cdot \frac{10^d B r}{R}.
\end{eqnarray*}
This gives that $k \le 2R/10^d r$, and therefore that the number of bad $w\in W'''$ is
$ \ge |W'''|-2R/10^dr > \frac{1}{2}|W'''| \ge R/10^d r$.

Let $W$ be the set of bad $w\in W'''$. Then $|W|\ge \frac{1}{2}|W'''|\ge \frac{1}{2\cdot 5^d}|W''|$.
Since $W''$ is connected, it follows that $W$ is connectible, so it satisfies the claim of the lemma.
\end{proof}

In the next subsection we prove the following theorem.

\begin{thm} \label{probbad} There exist constants $C,c>0$ such that
for any family $W \subseteq S$ we have
\begin{equation} \label{eq:probfamilybad}
\p\bigg(
\Omega_3^c \cap \Big\{ \text{all }w\in W\text{ are bad}\Big\}
\bigg) \le C e^{-c |W| \log R}.
\end{equation}
\end{thm}

Before turning to the proof of Theorem \ref{probbad}, here's how to prove Theorem \ref{thm3} using it.

\begin{proof}[Proof of Theorem \ref{thm3}]
Let $\Sigma$ be the set of connectible families $W\subseteq S$ with $|W|\ge R/10^d r$. By
Lemmas \ref{fewconnectible}, \ref{exceptionalomega} and \ref{manybad} we have that
\begin{eqnarray*}
\p(E_R) &\le& \p(\Omega_1) + \p(\Omega_2) + \p(\Omega_3) + \p(E_R \cap \Omega_1^c\cap\Omega_2^c \cap\Omega_3^c) \\ &\le&
Ce^{-c R^{1+\frac{1}{100d^2}}} + \sum_{W\in \Sigma}
 \p\bigg(\Omega_3^c \cap \Big\{ \text{all }w\in W\text{ are bad}\Big\}\bigg) \\ &\le&
 Ce^{-c R^{1+\frac{1}{100d^2}}} + \sum_{L\ge R/10^d r} \sum_{W\in \Sigma,\ |W|=L}
 \p\bigg(\Omega_3^c \cap \Big\{ \text{all }w\in W\text{ are bad}\Big\}\bigg) \\ &\le &
 Ce^{-c R^{1+\frac{1}{100d^2}}} + \sum_{L\ge R/10^d r} R^d C_{15}^L e^{-c L \log R} \\ &= &
O\bigg( e^{-c R (\log R)^{1-2/d}} \bigg).
\end{eqnarray*}
\end{proof}

\subsection{Bounding the probability for a collection of points to be bad \label{subsectionprobbad}}

Fix a family $W \subseteq S$. Our goal is to prove the inequality \eqref{eq:probfamilybad}. First, note that we may assume without loss of generality that the set $W$ is $12r$-separated in the infinity-norm;
that is, that for any $w,w'\in W$ with $w\neq w'$ we have that
$||w-w'||_\infty \ge 12 r$. Otherwise, replace $W$ with a
$12r$-separated
subset of it which has cardinality $\ge|W|/25^d$ (as in the proof of Lemma \ref{manybad} above),
and prove \eqref{eq:probfamilybad} for that subset. Throughout this subsection, we assume $W\subset S$ is a $12r$-separated family.

The next lemma can be deduced easily from the Besicovich covering lemma (see \cite{mattila}). For completeness we include a short direct proof.

\begin{lem} \label{geometric}
Given a set of $N$ balls $(B(x_i,r_i))_{i=1}^N$, where $r_i>1$ for all $i$ and $|x_i-x_j|>1$ for all $i,j$,
there exists a subset $(B(x_{i_j},r_{i_j}))_{j=1}^m$ of pairwise-disjoint balls satisfying
\begin{equation} \label{eq:geomet}
  \sum_{j=1}^m r_{i_j}^d \ge 6^{-d} N.
\end{equation}
\end{lem}

\begin{proof}
Assume that the radii $r_i$ are arranged in decreasing order. Construct the subsequence $(x_{i_j})_j$ sequentially as follows: $i_1 = 1$, and if we defined $x_{i_1},x_{i_2},\ldots,x_{i_t}$, take $i_{t+1}$ to be the least index $i > i_t$ such that the ball $B(x_i,r_i)$ is disjoint from $\cup_{j=1}^t B(x_{i_j},r_{i_j})$, or, if there is no such $i$, set $m=t$ and terminate. In the last step $t=m$, because the radii are decreasing, the fact that there was no index $i$ satisfying the requirements implies that the set $\cup_{j=1}^m B(x_{i_j},2r_{i_j})$ contains
all the points $x_i$. Therefore, since $r_i>1$ for all $i$, we have that
$$ \bigcup_{i=1}^N B(x_i,1/2) \subset \bigcup_{j=1}^m B(x_{i_j},3r_{i_j}). $$
The balls $\big(B(x_i,1/2)\big)_{i=1}^N$ are pairwise disjoint, so comparing the volumes of
both sides we get
$$ \kappa_d 2^{-d} N \le \kappa_d \sum_{j=1}^m 3^d r_{i_j}^d, $$
which finishes the proof.
\end{proof}

For each $w \in S$, define an event
\begin{eqnarray*}
\Omega_{4,w} &=&
\Big\{ \max_{x\in Q(w,3r)} \Big| U\big(x \big| B(x,3R^{1/d})\setminus B(x,3r)
 \big) \Big| \ge \frac{B}{4} \Big\}  \\
& &  \cup
  \Big\{ \max_{x\in Q(w,3r)}
    \Big| D_1 F\big(x \big| B(x,3R^{1/d})\setminus B(x,3r) \big) \Big| \ge \frac{\varepsilon}{4 \sqrt{d} \rho}\Big\}.
\end{eqnarray*}
Define a random set (depending on the fixed family $W$)
$$ {\cal W}_0 = \Big\{ w\in W : \Omega_{4,w}\text{ occurred} \Big\}.$$
Define an event (again depending on $W$)
$$ \Omega_{5,W} = \Big\{ |{\cal W}_0| > \frac{1}{2}|W| \Big\}.$$

\begin{lem} \label{dilutiontheorem}
Denote $\delta=\frac{1}{50d}$. For some constants $C,c>0$ depending only on $d$ we have that
$$  \p\big( \Omega_{5,W}\big) \le C e^{-c |W| R^\delta}
$$
\end{lem}

\begin{proof}
$$ \p(\Omega_{5,W}) \qquad\qquad\qquad\qquad\qquad\qquad\qquad\qquad\qquad\qquad\qquad\qquad\qquad
\qquad\qquad\qquad$$

\vspace{-20.0pt}
\begin{eqnarray*}
&=& P
\Big( \exists \textrm{ subcollection }W' \subset W\textrm{ with }|W'|> |W|/2
\textrm{ and }\bigcap_{w \in W'} \Omega_{4,w}\textrm{ holds }\Big) \\
&\le& 2^{|W|} \max_{W'\subset W,\ |W'|>|W|/2} \p\bigg(\bigcap_{w \in W'} \Omega_{4,w}\bigg).
\end{eqnarray*}
Therefore it's enough to prove that for some constants $C,c>0$,
for any subcollection $W'\subset W$ we have
\begin{equation} \label{eq:enough}
\p\bigg(\bigcap_{w \in W'} \Omega_{4,w}\bigg) \le C e^{-c|W'| R^\delta}.
\end{equation}
Fix a subcollection $W' \subset W$. Denote $\alpha = \frac{1}{20(d^2+1)}$.
Define a finite sequence of scales
\begin{equation*}
R^{1/d} = l_1 > l_2 > l_3 > \cdots > l_K = r
\end{equation*}
where we do not care about the precise values of the $l_i$ and only require that for each $i$ we have
\begin{equation} \label{eq:scale_relation}
l_{i+1} \ge \frac{l_i}{R^\alpha}
\end{equation}
and that $K$ is a constant depending only on $d$; for example, it is possible to define
such $l_i$ with $K=40d$.

For each $w\in W$ and each $1\le i<K$ define the event
\begin{eqnarray*}
A_w^i &=&
\Big\{ \max_{x\in Q(w,3r)}
\Big| U\big(x \big| B(x,3l_i)\setminus B(x,3l_{i+1})
 \big) \Big| \ge \frac{B}{4(K-1)} \Big\}  \\
& &  \cup
  \Big\{ \max_{x\in Q(w,3r)}
    \Big| D_1 F\big(x \big| B(x,3l_i)\setminus B(x,3l_{i+1}) \big) \Big| \ge \frac{\varepsilon}{4(K-1)\sqrt{d} \rho}\Big\}.
\end{eqnarray*}
We have, using Theorem \ref{largedeviations3}, that for some constants $c_2,c_3>0$
the estimate
\begin{equation*}
\P(A_w^i) \le Cr^d \exp\bigg(-c_1l_{i+1}^{d-2}B\log \left(\frac{c_2 B}{l_{i+1}^2}\right)\bigg) +
Cr^d \exp\left(-c_1l_{i+1}^d\frac{\eps}{\rho}\log(\frac{c_3 \eps}{\rho})\right)
\end{equation*}
holds. Using  \eqref{eq:scale_relation} and substituting the values $B=R^{8/9}$, $r=\Delta\cdot (\log R)^{2/d}$ and $\frac{\eps}{\rho} = R^{d/(10(d^2+1))}\log R$ we obtain
\begin{eqnarray} \label{eq:annulus_exceptional_event_estimate}
\P(A_w^i) &\le& C\exp\left(-c_4 l_i^{d-2}R^{8/9-\alpha(d-2)}\right)
+ C\exp\left(-c_4 l_i^d R^{d/10(d^2+1)-d\alpha}\right) \nonumber\\ &\le&
C\exp\left(-c_4 l_i^d R^{8/9 - \alpha(d-2)-2/d}\right) + C\exp(-c_4l_i^d R^\delta) \nonumber\\ &\le&
C\exp\left(-c_4 l_i^d R^{\delta}\right).
\end{eqnarray}
Note also that  $\Omega_{4,w} \subseteq \bigcup_{i=1}^{K-1} A_w^i$.
Therefore
\begin{eqnarray} \nonumber
\p\left( \bigcap_{w\in W'} \Omega_{4,w} \right) &\le& \p\left( \bigcap_{w\in W'} \bigcup_{i=1}^{K-1}
A_w^i \right)
= \p\left( \bigcup_{i:W'\to \{1,2,\ldots,K-1\}} \bigcap_{w\in W'} A_w^{i(w)} \right)
\\ &\le& K^{|W'|} \max_{i:W'\to \{1,2,\ldots,K-1\}} \p\left(\bigcap_{w\in W'} A_w^{i(w)} \right).
\label{eq:continuing}
\end{eqnarray}
Fix a function $i:W'\to \{1,2,\ldots,K-1\}$. We extract from the family of events $\left\{ A_w^{i(w)} \right\}_{w\in W'}$ a subfamily $\left\{ A_w^{i(w)} \right\}_{w\in W''}$ of \emph{independent} events, by using Lemma \ref{geometric}. By the definition of $A_w^i$, such a subfamily will be independent if the balls $\Big(B(w,7 \sqrt{d}l_{i(w)})\Big)_{w\in W''}$
are disjoint. By Lemma \ref{geometric} we can obtain such a subfamily with $\sum_{w\in W''} l_{i(w)}^d \ge (42\sqrt{d})^{-d} |W'|$. This gives, continuing \eqref{eq:continuing} and using \eqref{eq:annulus_exceptional_event_estimate},
that
\begin{eqnarray*}
\p\left( \bigcap_{w\in W'} \Omega_{4,w} \right) &\le& K^{|W'|}C^{|W'|} \max_{i:W'\to \{1,2,\ldots,K-1\}}
 C \exp\bigg( -c_0 \sum_{w\in W''} l_{i(w)}^d R^\delta \bigg)
\\ &\le&  C \exp\bigg( -c_0 \cdot (42\sqrt{d})^{-d} \cdot |W'| R^\delta \bigg).
 \end{eqnarray*}
 This proves \eqref{eq:enough} and finishes the proof of the lemma.
\end{proof}

For each $w\in S$ introduce two subgrids of points
\begin{eqnarray*}
S_w &=& s\Z \cap (Q(w,2r)\setminus Q(w,r)), \\
T_w &=& \rho\Z\cap (Q(w,2r)\setminus Q(w,r)).
\end{eqnarray*}

For $w\in S$, two subgrid points $x,x'\in S_w$ are called adjacent if $|x-x'|=s$. A {\bf chain} is a sequence of points such that each two consecutive points are adjacent. A point $x\in S_w$ is called an {\bf inner point of $S_w$} if $\textrm{dist}(x,\partial Q(w,r)) < s$, and it is called an {\bf outer point of $S_w$} if $\textrm{dist}(x,\partial Q(w,2r)) < s$.

\begin{defi} A point $x\in \R^d$ is called {\bf $\alpha$-crowded} if $Q(x,\alpha s)$ contains a star. A point $w\in S$ is called {\bf percolating} if there exists a chain of distinct points $(x_i)_{0\le i\le k} \subset S_w$ such that
$x_0$ is an inner point of $S_w$, $x_k$ is an outer point of $S_w$, and at least a $9/10$-fraction of the $x_i$'s are $8$-crowded.
\end{defi}

With $W \subset S$ as above a $12r$-separated family, define a random
set (again depending on $W$)
$$ {\cal W}_1 = \Big\{ w\in W : w\text{ is percolating} \Big\}.
$$
Define an event
$$ \Omega_{6,W} = \Big\{ |{\cal W}_1| > \frac{1}{2}|W| \Big\}.$$

\begin{lem} \label{seconddilution}
For some constants $C,c,\alpha >0$ depending only on $d$ we have
that
$$ \p\big(\Omega_{6,W}\big) \le C e^{-c |W| R^\alpha}. $$
\end{lem}

\begin{proof} For any $x\in\R^d$ we have for some constant $c_1>0$ that
$$\p(x\textrm{ $8$-crowded}) = 1-e^{-(16s)^d} \le c_1 s^d. $$
Fix a $w\in S$. For each $k$, the number of chains of distinct points
$(x_i)_{0\le i\le k} \subset S_w$ such that
$x_0$ is an inner point of $S_w$ and $x_k$ is an outer point of $S_w$
is at most $2d\cdot (2r/s)^{d-1}\cdot (2d)^k$.
Note that such chains can only exist if $k \ge r/s$,
so in particular, the number of such chains
 is $\le C_2^k$ for some constant $C_2>0$. For each such chain ${\cal
C}$, there are $\le 2^k$ subsets of it of size at least $9k/10$. Given
such a chain ${\cal C}$ and such a subset ${\cal C}'$ of it,
one may choose using a greedy method (as in Lemma \ref{manybad}
a further subset $(x_i')_{0\le i\le k'}$ of ${\cal C}'$ with
$k' \ge 9k/(10\cdot 33^d)$, such that for each $i \neq j$ we have that
$||x_i'-x_j'||_\infty\ge 16 s$.
Therefore the events
 $\big(\{x_i'\textrm{ is $8$-crowded}\}\big)_{0\le i\le k'}$ are independent. Therefore we have the bound
\begin{eqnarray*}
\p(w\textrm{ percolating}) &\le& \sum_{k=\lfloor \frac{r}{s}\rfloor
 }^\infty 2^k
 C_2^k \left(c_1 s^{d}\right)^{\frac{9}{10\cdot 33^d}\cdot k}
\le C_3 \left( C_4 s^d \right)^{\left(\frac{C_5 r}{s} \right)} \le C_6 e^{-c_7 R^\alpha}
\end{eqnarray*}
for some $C_3, C_4,C_5,C_6, c_7, \alpha>0$ and all $R$ large enough.
Now, because of the assumptions that the points of $W$ are $12r$-separated, the events $\big( \{w\textrm{ percolating}\}\big)_{w\in W}$ are independent. Therefore
\begin{eqnarray*}
\p\Big(\Omega_{6,W}\Big)
& = &
\p\bigg( \exists\textrm{ subcollection }W'\subset W\textrm{ with }|W'|> |W|/2 \\ & & \qquad
\textrm{ and all }w\in W'\textrm{ are percolating }\bigg) \\ &\le &
2^{|W|} \left(C_6 e^{-c_7 R^\alpha} \right)^{|W|/2} \le C_8 e^{-c_9 |W| R^{\alpha}}.
\end{eqnarray*}
\end{proof}

\begin{defi} If $w\in S$, a point $x \in T_w$ is called {\bf black} if there exists a point
$y \in Q(x,2 \rho)$ such that either $|F(y)| < \frac{10000^d B}{R}$ or
$ U\big(y\ \big|\ B(y,3R^{1/d})\big) < -\frac{B}{2}$.
\end{defi}

\begin{lem} \label{haveblack} 
If $w\in S$ is bad and not percolating, then
the subgrid $T_w$ contains at least $\frac{r}{100^d \rho}$ points which are black and not $4$-crowded.
\end{lem}

\begin{proof} Let $w\in S$ be bad and not percolating. Let $\gamma$ be a gravitational flow curve connecting $\partial Q(w,r)$ with $\partial Q(w,2r)$ such that at least one of the conditions (1), (2) in Definition \ref{defibad} holds. Let ${\cal I}$ be the set of points $x\in T_w$ such that some segment of $\gamma$  crosses from $\partial Q(x,\rho)$ to $\partial Q(x,2\rho)$. Let ${\cal J}$ be the set of points $x\in S_w$ such that some segment of $\gamma$ crosses from $\partial Q(x,s)$ to $\partial Q(x,2s)$. Note that
$|{\cal J}| \ge \frac{r}{s}-2\ge \frac{0.95r}{s}$ (with the rightmost inequality holding for $R$ sufficiently large), and for each $x\in {\cal J}$ we have that $|{\cal I} \cap (Q(x,2s)\setminus Q(x,s))| \ge \frac{s}{\rho}-2 \ge \frac{0.9s}{\rho}$, with the rightmost inequality holding for $R$ sufficiently large.

First, we prove that the number $\nu$ of $x\in {\cal I}$ which are not black is at most $\frac{r}{100\cdot 11^d \rho}$. If condition (1) in Definition \ref{defibad} holds for $\gamma$, then for all $x\in {\cal I}$, for some $z\in Q(x,2\rho)\setminus Q(x,\rho)$ which is in the range of $\gamma$ we have that $U\big(z\ \big|\ B(z,3R^{1/d})\big)<-B/2$, so $x$ is black. So we have shown that if condition (1) in Definition \ref{defibad} holds,  all $x\in {\cal I}$ are black. The other possibility is that condition (2) in Definition \ref{defibad} holds for $\gamma$. In that case, denote by $\gamma'$ the union of those segments of $\gamma$ crossing from $\partial Q(x,\rho)$ to $\partial Q(x,2\rho)$ for those $x\in {\cal I}$ which are not black.
It is not difficult to see that $\text{len}(\gamma')\ge \rho
\frac{\nu}{13^d}$ (this uses a similar argument to the ones used at
the beginning of this subsection and in the proofs of Lemmas \ref{manybad},
\ref{seconddilution}).
So, because of the definition of blackness we have that
\begin{eqnarray*}
\frac{10^d Br}{R}&>& \int_\gamma |F(x)|\cdot |dx| \ge \int_{\gamma'} |F(x)|\cdot |dx| \ge \frac{10000^d B}{R}\textrm{len}(\gamma')
\\ &\ge&\frac{10000^d B}{R}\cdot \frac{\rho\cdot \nu}{13^d},
\end{eqnarray*}
and therefore $\nu \le \frac{r}{100\cdot 11^d \rho}$, as claimed.

Next, because of the assumption that $w$ is not percolating, at least
a $1/10$-fraction of the points $x\in {\cal J}$ are not $8$-crowded
(minus $2$ to account for boundary effects),
or in absolute terms at least $\frac{0.94r}{10s}$ points in ${\cal J}$.
As in the proof of Lemma \ref{seconddilution} above, we can choose a further subset ${\cal J}' \subset {\cal J}$ of those points of ${\cal J}$ which are not $8$-crowded which is $5s$-separated and such that
$|{\cal J}'| \ge \frac{0.94r}{10\cdot 11^d s}$.
By the remark in the first paragraph of the proof, for each such $x\in {\cal J}'$ we have at least $\frac{0.9s}{\rho}$ points
$y \in {\cal I}\cap (Q(x,2s)\setminus Q(x,s))$, and these $y$ are not $4$-crowded. That gives a total of at least $\frac{0.8r}{10\cdot 11^d\rho}$ points $y\in {\cal I}$ which are not $4$-crowded, and these points
are all distinct because ${\cal J}'$ is $5s$-separated. Since as we proved above at most
$\frac{r}{100\cdot 11^d\rho}$ of them are not black, it follows that there are at least
$\frac{r}{100^d\rho}$ points $y\in {\cal I}$ which are black and not $4$-crowded, as claimed.
\end{proof}

Let $W \subset S$ be a $12r$-separated family. Denote
\begin{multline*}
\Omega_{7,W} = \Omega_3^c \cap \bigcap_{w\in W} \Omega_{4,w}^c  \\ \cap \bigcap_{w\in W}
\Big\{ \textrm{at least }\frac{r}{100^d \rho}\textrm{ points }x\in T_w\textrm{ are black and not 4-crowded}
\Big\}.
\end{multline*}

\begin{lem} \label{pwbound} For some constants $C,c>0$, we have
$$ \p(\Omega_{7,W}) \le C e^{-c |W| \log R }. $$
\end{lem}

\begin{proof}
Let $N=|W|\le |S|\le (4R/r)^d$,
and let $W=\{w_1,w_2,\ldots,w_N\}$ be some arbitrary ordering of the points
of $W$, and define a random variable
$$X_{W} = \#\Big\{ (x_1,x_2,\ldots,x_N) \in \prod_{j=1}^N T_{w_j}\ \Big|\ x_j\textrm{ are all
black and not 4-crowded} \Big\}. $$
Then, by Markov's inequality,\\
\vbox{
\begin{eqnarray} \nonumber
\p(\Omega_{7,W}) &\le & \p\bigg( \bigg\{X_{W} \ge \left(\frac{r}{100^d \rho} \right)^N \bigg\} \cap
\Omega_3^c \cap \bigcap_{w\in W} \Omega_{4,w}^c \bigg) 
\qquad\qquad\qquad
\end{eqnarray}

\vspace{-20.0pt}
\begin{eqnarray} \nonumber
&\le& \left(\frac{r}{100^d \rho}\right)^{-N}  \!\!\!\!\!\!\! \sum_{(x_1,\ldots,x_N)\in\prod_{j=1}^N T_{w_j}}
\!\!\!\!\!\! \p\bigg( \bigcap_{j=1}^N \Big( \Omega_3^c\cap
\Omega_{4,{w_j}}^c \cap \{x_j\textrm{\ black, not 4-crowded}\} \Big)
\bigg) \\&=&
\left(\frac{r}{100^d \rho}\right)^{-N}  \!\!\!\!\!\!\! \sum_{(x_1,\ldots,x_N)\in\prod_{j=1}^N T_{w_j}}
\E\left[ \prod_{j=1}^N 1_{A_{j,x_j}}\right], \label{eq:pw}
\end{eqnarray}
}
where $A_{j,x_j} = \Omega_3^c\cap \Omega_{4,{w_j}}^c \cap \{x_j\textrm{ is black and not 4-crowded}\}$.
Denote $\lambda = 2(d/\kappa_d)^{1/d}(\log R)^{1/d}$. Fixing $(x_1,\ldots,x_N)$ for the moment to simplify the notation, for $1\le j\le N$ denote
\begin{eqnarray*}
A_j' &=& \bigg\{ \max_{x\in Q(x_j,2\rho)}  \Big| U\big( x\ \big|\ B(x,3r)\setminus B(x,3s) \big)
\Big| > \frac{B}{4} \bigg\} \\
  & & \cup
              \bigg\{ \max_{x\in Q(x_j,2\rho)} \Big| D_1 F\big(x\ \big|\ B(x,3r)\setminus B(x,3s) \big)
              \Big| > \frac{\varepsilon}{16\sqrt{d} \rho} \bigg\} \\
 & & \cup \bigg\{\textrm{ no star in }B(x_j,\lambda) \bigg\}, \\
A_j'' &=&  \bigg\{ |F(x_j)|<\varepsilon,\textrm{ some stars in }B(x_j,\lambda)\bigg\}.
\end{eqnarray*}
We claim that
\begin{equation}\label{eq:ajcontained}
A_{j,x_j} \subset A_j' \cup A_j''.
\end{equation}
Here is the proof: Assume $A_{j,x_j}$ holds. If for some $y\in Q(x_j,2 \rho)$ we have that $|F(y)|<\frac{10000^d B}{R}$,
which for $R$ sufficiently large implies $|F(y)|<\frac{\varepsilon}{16}$, then: Either $A_j'$ occurred, or, if not, then from the definitions of $4$-crowdedness and of the events
$\Omega_3^c, \Omega_{4,w_j}^c, (A_j')^c$ it follows that for $R$ large enough,
\begin{eqnarray*}
\max_{x\in Q(x_j,2\rho)} |D_1 F(x)| &\le&
\max_{x\in Q(x_j,2\rho)} \Big|D_1 F(x\ |\ B(x,3s))\Big|
\\ & & +
\max_{x\in Q(x_j,2\rho)} \Big|D_1 F(x\ |\ B(x,3r)\setminus B(x,3s))\Big|
\\ & & +
\max_{x\in Q(x_j,2\rho)} \Big|D_1 F(x\ |\ B(x,3R^{1/d})\setminus B(x,3r))\Big|
\\ & &+
\max_{x\in Q(x_j,2\rho)} \Big|D_1 F(x\ |\ \R^d\setminus B(x,3R^{1/d}))\Big|
\\ & \le& \sqrt{d}\kappa_d + \frac{\eps}{16\sqrt{d}\rho} +
\frac{\eps}{4\sqrt{d}\rho} +
\frac{\eps}{8\sqrt{d}\rho} < \frac{15\eps}{32\sqrt{d}\rho }
\end{eqnarray*}
(the term $\sqrt{d}\kappa_d$ comes from eq. \eqref{eq:numberedeq2}, note that it is $<< \eps/\rho$ for large $R$).
Therefore $|F(x_j)|<\frac{\eps}{16}+2\rho\sqrt{d}
\cdot \frac{15\eps}{32\sqrt{d}\rho} =\eps$, and since $A_j'$ did not occur there are stars in $B(x_j,\lambda)$ and therefore $A_j''$ occurred.

The other possibility by blackness of $x_j$ is that for some $y\in Q(x_j,2\rho)$ we have that
$U\Big(y\ \Big|\ B(y,3R^{1/d})\Big) < -\frac{B}{2}$. Then, because of $\Omega_{4,w_j}^c$ we also have that $U\Big(y\ \Big|\ B(y,3r)\Big) < -\frac{B}{4}$, and because $x_j$ is not 4-crowded we can write equivalently
$U\Big(y\ \Big|\ B(y,3r)\setminus B(y,3s)\Big) < -\frac{B}{4}$, so $A_j'$ occurred. This completes the proof of \eqref{eq:ajcontained}.

Note that each $A_j'$ is measurable with respect to the locations of the
stars in $Q(w_j,5.5r)$. Therefore, because the $w_j$ are $12r$-separated,
$(A_j')_j$ is an independent family of events. By Theorem \ref{largedeviations3} we have
for each $j$ that
\begin{eqnarray*}
\p(A_j')& \le& Ce^{-c s^{d-2}B\log\left(\frac{c_3 B}{s^2}\right)}
+ C e^{-c s^d \frac{\varepsilon}{\rho}\log \left(\frac{c_3 \varepsilon}{\rho}
\right) }
+ e^{-\kappa_d \lambda^d} \\ & \le &
 \frac{1}{R^d}\ \ \ \textrm{(for }R\text{ sufficiently large.)}
 \end{eqnarray*}
Also, for any integers $1\le b_1 < b_2< \ldots < b_i \le N$, by Theorem \ref{joint_force_density_thm}
we have that almost surely
$$ \p\Big(A_{b_1}''\cap \ldots \cap A_{b_i}''\ \Big|\ {\cal F}_{b_1,\ldots,b_i}
 \Big) \le \left(\varepsilon^d C \lambda^{d^2-d}\right)^i, $$
where ${\cal F}_{b_1,\ldots,b_i}$ is the $\sigma$-algebra generated by the
locations of the stars in $\Big(B(x_{b_1},\lambda)\cup B(x_{b_2},\lambda)\cup \ldots
\cup B(x_{b_i},\lambda)\Big)^c$, provided that the assumptions of that theorem hold; because of the values chosen for the parameters $r$ and $\lambda$, this is true if $\Delta$ is chosen to be a sufficiently large constant. This is the only place where the value of $\Delta$ is important. So we can write, for sufficiently large $R$,
$$ \E\left[ \prod_{j=1}^N 1_{A_{j,x_j}}\right] \le
 \E\left[ \prod_{j=1}^N \left(1_{A_j'}+1_{A_j''} \right) \right]  \qquad\qquad\qquad\qquad
 \qquad\qquad\qquad\qquad$$

\vspace{-10.0pt}
\begin{eqnarray} &\le&
 \sum_{i=0}^N \binom{N}{i}
  \max_{\scriptsize \begin{array}{l}  a_1< \ldots< a_{N-i} \nonumber \\b_1<\ldots<b_i
\\  \forall j,\ell\ \ a_j \neq b_\ell \end{array}}
\p\Big( A_{a_1}'\cap \ldots \cap A_{a_{N-i}}' \bigcap
              A_{b_1}''\cap \ldots \cap A_{b_i}'' \Big) \nonumber \\
&\le&
 \sum_{i=0}^N \binom{N}{i} \frac{1}{R^{d (N-i)}}
  \max_{\scriptsize \begin{array}{l}  a_1< \ldots< a_{N-i} \nonumber \\b_1<\ldots<b_i
  \\  \forall j,\ell\ \ a_j \neq b_\ell\end{array}}  \!\!\!\!\!
\p\Big( A_{b_1}''\cap \ldots \cap A_{b_i}''\ \Big|\
A_{a_1}'\cap \ldots \cap A_{a_{N-i}}'
\Big) \nonumber \\
 &\le&
\sum_{i=0}^N \binom{N}{i} \frac{1}{R^{d (N-i)}}
\left(\varepsilon^d C \lambda^{d^2-d}\right)^i \nonumber\\
&\le& \left(\frac{1}{R^d} + R^{-\frac{d}{10}+\frac{d^2}{10(d^2+1)}+\frac{1}{20(d^2+1)}}\right)^{N}
\le \left(2 R^{-\frac{d}{10}+\frac{d^2}{10(d^2+1)}+\frac{1}{20(d^2+1)}} \right)^N \label{eq:thiseq}
\end{eqnarray}
(since for large $R$ the polylogarithmic factor in $\eps^d C \lambda^{d^2-d}$ can be bounded by
$R^{\frac{1}{20(d^2+1)}}$).
Now \eqref{eq:thiseq} gives, using \eqref{eq:pw}, that for sufficiently large $R$ we have
\begin{eqnarray*}
\p(\Omega_{7,W}) &\le& \left(\frac{r}{100^d\rho}\right)^{-N} \left(\frac{r}{\rho}\right)^{dN}
\left( 2 R^{-\frac{d}{10}+\frac{d^2}{10(d^2+1)}+\frac{1}{20(d^2+1)}} \right)^N \\ &\le &
\bigg[Cr^{d-1} R^{\frac{d-1}{10}}\left(2 R^{-\frac{d}{10}+\frac{d^2}{10(d^2+1)}+\frac{1}{20(d^2+1)}}\right)\bigg]^N \\ &\le&
C e^{-\frac{1}{40(d^2+1)} N \log R} \le C e^{-c |W| \log R}.
\end{eqnarray*}
\end{proof}

The only step remaining to complete the proof of Theorem \ref{thm3} in dimensions 5 and higher is the following.

\begin{proof}[Proof of Theorem \ref{probbad}.]
As noted above, we assume that $W$ is a $12r$-separated family. For an event $A$ denote $\p'(A)=\p(A\cap 
\Omega_3^c)$. The idea of the proof is roughly as follows. Because of Lemma \ref{dilutiontheorem}, we may replace $W$ by a subset $W'$ so that $|W'|\ge |W|/2$ and such that $\cap_{w\in W'} \Omega_{4,w}^c$ occurred. Because of Lemma \ref{seconddilution}, we may replace $W'$ by a further subset $W'' \subset W'$ so that $|W''| \ge |W'|/2$ and all $w\in W''$ are not percolating.
Finally, by Lemmas \ref{haveblack} and \ref{pwbound}, the probability that all $w\in W''$ are bad is $\le C e^{-c |W| \log R}$.

Formally, we have
\begin{eqnarray*}
\p'\left(\bigcap_{w\in W} \Big\{w\text{ bad}\Big\}\right) &\le& \p'(\Omega_{5,W}) +
\p'\left(\Omega_{5,W}^c \cap \bigcap_{w\in W} \Big\{w\text{ bad}\Big\}\right) 
\ \ \
\end{eqnarray*}

\vspace{-10.0pt}
\begin{eqnarray} \nonumber \ \ &\le&
Ce^{-c|W| R^\delta} + \sum_{\scriptsize \begin{array}{c}W'\subset W, \\|W'|\ge |W|/2\end{array}}
\p'\left( \bigcap_{w\in W'} \Big(\Big\{w\text{ bad}\Big\}\cap \Omega_{4,w}^c \Big) \right) \\ &\le &
Ce^{-c|W| R^\delta} + 2^{|W|} \max_{\scriptsize \begin{array}{c}W'\subset W,\\ |W'|\ge |W|/2\end{array}}
\p'\left( \bigcap_{w\in W'} \Big(\Big\{w\text{ bad}\Big\}\cap \Omega_{4,w}^c \Big) \right).
\label{eq:combining1}
\end{eqnarray}
For any $W'$ we have
\begin{eqnarray*}
\p'\left( \bigcap_{w\in W'} \Big(\Big\{w\text{ bad}\Big\}\cap \Omega_{4,w}^c \Big) \right)
\qquad\qquad\qquad\qquad\qquad\qquad\qquad\qquad
\end{eqnarray*}

\vspace{-15.0pt}
\begin{eqnarray} \ \ &\le & \nonumber
\p'\left(\Omega_{6,W'}\right) + \p'\left(\Omega_{6,W'}^c \cap
\bigcap_{w\in W'} \Big(\Big\{w\text{ bad}\Big\}\cap \Omega_{4,w}^c \Big) \right) \\ &\le &
Ce^{-c|W'|R^\alpha}  \nonumber \\ & &+
\sum_{\scriptsize \begin{array}{c}W''\subset W', \\|W''|\ge |W'|/2\end{array}}
\p'\left( \bigcap_{w\in W''} \Big(\Big\{w\text{ bad, not percolating}\Big\}\cap \Omega_{4,w}^c \Big) \right)
\nonumber \\ &\le &
Ce^{-c|W'|R^\alpha}  \nonumber \\ & &+ 2^{|W'|}
\!\!\!\!\!\!\!\!\!\max_{\scriptsize \begin{array}{c}W''\subset W',  \\|W''|\ge |W'|/2\end{array}}
\!\!\!\!\!\!
\p'\left( \bigcap_{w\in W''} \Big(\Big\{w\text{ bad, not percolating}\Big\}\cap \Omega_{4,w}^c \Big) \right)\!.
\label{eq:combining2}
\end{eqnarray}
For any $W''$ we have, by Lemmas \ref{haveblack}, \ref{pwbound},
\begin{eqnarray} \label{eq:combining3}
\p'\left( \bigcap_{w\in W''} \Big(\Big\{w\text{ bad, not percolating}\Big\}\cap \Omega_{4,w}^c \Big) \right)
&\le& Ce^{-c |W''| \log R}.
\end{eqnarray}
Combining \eqref{eq:combining1}, \eqref{eq:combining2}, \eqref{eq:combining3}, and remembering that
$|W''| \ge |W'|/2 \ge |W|/4$, we get
$$ \p'\left(\bigcap_{w\in W} \Big\{w\text{ bad}\Big\}\right) \le C e^{-c |W| \log R}. $$
This completes the proof.
\end{proof}

\section{Proof of Theorem \ref{thm3} in dimensions 3 and 4 \label{sectiondim34}}

In this section, we prove Theorem \ref{thm3} in dimensions 3 and 4. Much of the proof in
dimensions 5 and higher remains unchanged. However, new large deviations estimates are
required, as well as the introduction of a new function, the {\bf potential difference function}.

\subsection{The potential difference function}

First, it is instructive to understand why the proof in Section \ref{sectionproofthm3} fails in dimensions 3 and 4. The difficulty is that the stationary potential function $U(x)$ cannot be defined. This can be seen from equation \eqref{eq:varianceupq}: in dimension $3$ the variance of $U_{p,q}$ diverges like a constant times $p$ as $p\to\infty$, and in dimension $4$ like a constant times $\log p$.

However, the proof in Section \ref{sectionproofthm3} for the most part does not use the full stationary potential. After discarding the atypical events $\Omega_1, \Omega_2$, only the partial potential $U(x\ |\ A)$ is used
for various sets $A\subset B(x,3R^{1/d})$. So, to adapt the proof to dimensions 3 and 4, we replace these events with suitable adaptations of them.

Assume for the rest of this section that $d=3$ or $4$. Define the potential difference function $\pd:\R^d\times \R^d\to \R$ by
\begin{equation}\label{eq:pd}
\pd(x,y) = \frac{1}{d-2}\sum_{|z_i|\uparrow} \left( \frac{-1}{|z_i-y|^{d-2}}- \frac{-1}{|z_i-x|^{d-2}}\right) +
\frac{\kappa_d}{2}\left(|x|^2-|y|^2\right)
\end{equation}
where the sum is in order of increasing $|z_i|$.
We need to check that this sum converges a.s. This is true because, defining
\begin{equation}\label{eq:deltau}
\pd_{p,q}(x,y) = \sum_{\scriptsize \begin{array}{c} |z_i|<p,\\ |z_i-x|>q, |z_i-y|>q\end{array}} \left( \frac{-1}{|z_i-y|^{d-2}}- \frac{-1}{|z_i-x|^{d-2}}\right),
\end{equation}
it is easy to check as in Section \ref{sectionpotential} that if $|x|,|y|<p-q$ and $|x-y|>2q$ then
\begin{eqnarray*}
\E[\pd_{p,q}(x,y)] &=& \int_{B(0,p)\setminus B(x,q)\setminus B(y,q)} \left( \frac{-1}{|z-y|^{d-2}}
- \frac{-1}{|z-x|^{d-2}} \right)dz \\ &=& \frac{(d-2)\kappa_d}{2}\left( |y|^2-|x|^2 \right),
\end{eqnarray*}

\vspace{-15.0pt}
\begin{eqnarray*}
\var\Big[\pd_{p,q}(x,y)-\pd_{p',q}(x,y)\Big] &=& O\left(\frac{1}{p^{d-2}}+\frac{1}{{p'}^{d-2}}\right), \quad p,p'\to\infty.
\end{eqnarray*}
Similarly, using the methods of Section \ref{sectionproofprop1}, it is not difficult to prove the following.

\begin{lem}
The series \eqref{eq:pd} converges a.s.\ simultaneously for all $x,y\in \R^d\setminus \{z_i\}_i$,
and defines a centered process that a.s.\ is differentiable where it is defined and satisfies
$ \nabla_{\!\!x}\,\pd(x,y) = F(x)$, $\nabla_{\!\!y}\,\pd(x,y) = -F(y)$.
\end{lem}

\paragraph{Remark.} The potential difference can in fact be defined for all dimensions $d\ge 3$, and for dimensions $d\ge 5$ we have that $\pd(x,y) = U(y) - U(x)$.

\medskip For a bounded set $A\subset \R^d$, denote
\begin{eqnarray*}
\pd(x,y\ |\ A) &=& \frac{1}{d-2}
\sum_{z_i \in A} \left( \frac{-1}{|z_i-y|^{d-2}}- \frac{-1}{|z_i-x|^{d-2}}\right) \\ & &
-  \frac{1}{d-2}\int_A \left( \frac{-1}{|z-y|^{d-2}}
- \frac{-1}{|z-x|^{d-2}} \right)dz.
\end{eqnarray*}
For a set $A\subset \R^d$ whose complement is bounded, denote
$$ \pd(x,y\ |\ A) = \pd(x,y) - \pd(x,y\ |\ \R^d\setminus A). $$
Again, it can be verified that if $A\subset \R^d$ is an annulus of the
form $B(v,p)\setminus B(v,q)$, where $v\in\R^d$ and $0\le q<p\le
\infty$, then
$$ \nabla_x \pd(x,y\ |\ A) = F(x\ |\ A), \qquad \nabla_y \pd(x,y\ |\ A) =
-F(y\ |\ A). $$

\subsection{Large deviations estimates in dimensions 3 and 4}

The large deviations estimates which we prove in this subsection will complement the estimates in Section \ref{sectionlarge}.

\begin{thm} \label{largedeviations4}
In dimension $d=4$, there exist constants $C_1, c_2, c_3>0$ such that
for all $x,y\in \R^d$ and $p>q>2$ satisfying $p>3q, |x|,|y|<p/2$, and $|x-y|>3q$, we have that
\begin{equation}\label{eq:lrg4}
 \p\left(  \Big| \pd\big(x,y\ \big|\ B(0,p) \setminus
(B(x,q)\cup B(y,q)) \big)\Big| > t \right) \le C_1 e^{-c_2 q^2 t \log\left(\frac{c_3 t}{q^2}\right)}
\end{equation}
for all $t$ above a threshold that depends on $q$ and $|x-y|$, as follows:
\begin{equation}
\begin{array}{ll} t\ge C_1 q^2\text{ and } \\ t\ge C_1 q^2 \log\left(\frac{C_1 t}{q^2}\right)
  \log\left(\frac{|x-y|}{q}\right). 
  \end{array}
  \label{eq:lrgassumption1}
\end{equation}
Similarly, in dimension $d=3$, there exist constants $C_1, c_2, c_3>0$ such that
for all $x,y\in \R^d$, $p>q>2$ and $t>0$ satisfying $p>3q, |x|,|y|<p/2$, and $2<q < t < \frac{1}{3}|x-y|$
 we have that
\begin{equation}\label{eq:lrg44}
 \p\left(  \Big| \pd\big(x,y\ \big|\ B(0,p) \setminus
(B(x,q)\cup B(y,q)) \big)\Big| > t \right) \le C_1 e^{-\frac{c_2 t^2}{|x-y|}}.
\end{equation}
\end{thm}

\begin{proof} The proof is modeled after the proof of Theorem \ref{largedeviations1}. Fix $x, y, p$ and $q$. Denote $B = B(0,p)\setminus B(x,q)\setminus B(y,q)$. Let $W$ be a uniform random point in $B$, let $N$ be the number of stars in $B$, and let
$$\pd = \sum_{i\,:\ z_i\in B} \left(\frac{-1}{|z_i-y|^{d-2}}- \frac{-1}{|z_i-x|^{d-2}}\right)
= \pd_{p,q}(x,y) $$
as defined in \eqref{eq:deltau}, so that
$ \pd\big(x,y\ \big|\ B(0,p) \setminus (B(x,q)\cup B(y,q)) \big) = \frac{1}{d-2}(\pd-\E[\pd])$.
Let $$ V = |W-x|^{-(d-2)} - |W-y|^{-(d-2)}. $$
Then for any $u>0$ we have
\begin{eqnarray} \nonumber
\p(|V|>u) &\le& \p\left(\frac{1}{|W-x|^{d-2}}>u\right) +
\p\left(\frac{1}{|W-y|^{d-2}}>u\right) \\ &\le&
\frac{2\kappa_d}{\vol(B) u^{d/(d-2)}}. \label{eq:tailbound}
\end{eqnarray}
Suppose now $d=4$.
Noting that $|V| \le 2/q^2$, for any integer $k\ge 3$ we have exactly as in \eqref{eq:moment}  that
for some constant $C_7>0$,
\begin{eqnarray*}
\E(|V|^k)  &=& \int_0^{2/q^2} ku^{k-1}\p(|V|>u)du \le \frac{C_7}{\vol(B) q^{2k-4}},
\end{eqnarray*}
Evidently, we need a better tail bound for $|V|$ to get anything useful for $k=2$. To that end, note that
\begin{eqnarray*}
|V| &=& \frac{\big| |W-x|^2-|W-y|^2\big|}{|W-x|^2\cdot |W-y|^2} \le \frac{|x-y|(|W-x|+|W-y|)}{|W-x|^2
\cdot |W-y|^2} \\ &=&|x-y|\left( \frac{1}{|W-x|\cdot |W-y|^2}+\frac{1}{|W-y|\cdot |W-x|^2} \right).
\end{eqnarray*}
Now, if $|W-x|$ and $|W-y|$ are both bigger than $(2|x-y|/u)^{1/3}$, then a simple verification using the above inequality shows that $|V|\le u$. Thus,
\begin{eqnarray}\nonumber
\p(|V|>u) &\le& \p\bigg(|W-x|\le \left(\frac{2|x-y|}{u}\right)^{1/3}\bigg) \\ & + &\p\bigg(|W-y|\le
\left(\frac{2|x-y|}{u}\right)^{1/3}\bigg) \le C_8
\frac{|x-y|^{4/3}}{\vol(B) u^{4/3}}. \label{eq:tailbound2}
\end{eqnarray}
Combining the bounds from \eqref{eq:tailbound} and \eqref{eq:tailbound2} and using the assumption
that $|x-y|>3q$, we get that
$$
\E|V|^2 \le \int_0^{|x-y|^{-2}} 2u \frac{C_8 |x-y|^{4/3}}{\vol(B) u^{4/3}} du + \int_{|x-y|^{-2}}^{2q^{-2}}
2u\frac{2\kappa_4}{\vol(B)u^2} du \le \frac{C_9 \log\left(\frac{|x-y|}{q}\right)}{\vol(B)}.
$$
Thus, we have for any $\theta\ge 0$ that
\begin{eqnarray*}
\E(e^{\theta V}) &\le& 1 + \theta\E(V)+\theta^2 \frac{C_9 \log\left(\frac{|x-y|}{q}\right)}{2\vol(B)}
+ \frac{C_7 q^4}{\vol(B)} \sum_{k=3}^\infty \frac{(\theta/q^2)^k}{k!}
\\&\le& 1+ \theta\E(V) + \theta^2 \frac{C_9 \log\left(\frac{|x-y|}{q}\right)}{2 \vol(B)} +
\frac{C_7 q^4}{\vol(B)} e^{\theta/q^2}.
\end{eqnarray*}
Now proceed as in the proof of Theorem \ref{largedeviations1}. Conditionally on $N$, the stars in $B$ are a vector of $N$ i.i.d points distributed uniformly in $B$, and therefore
$$ \E\Big[e^{\theta \pd}\ \big|\ N\Big]
\le \left(1+ \theta\E(V) + \theta^2 \frac{C_9 \log\left(\frac{|x-y|}{q}\right)}{2\vol(B)} +
\frac{C_7 q^4}{\vol(B)} e^{\theta/q^2}
\right)^{N}. $$
This implies as before that
$$
\E\left(e^{\theta (\pd-\E(\pd)}\right) \le \exp\bigg(
\frac{1}{2}C_9 \theta^2 \log\left(\frac{|x-y|}{q}\right) + C_7 q^4 e^{\theta/q^2}\bigg),
$$
whence, for any $t\ge 0$ and $\theta\ge 0$ we have
$$ \p\left( \pd - \E[\pd] \ge t \right) \le
\exp\bigg(-\theta t +
\frac{1}{2}C_9 \theta^2 \log\left(\frac{|x-y|}{q}\right) + C_7 q^4 e^{\theta/q^2}\bigg).
$$
Take $\theta = q^2 \log(t/C_7 q^2)$. Then, if we assume \eqref{eq:lrgassumption1} for some sufficiently large constant $C_1>0$, we get that
\begin{eqnarray*}
\frac{1}{2}C_9 \theta^2 \log\left(\frac{|x-y|}{q}\right)  
\le \frac{\theta t}{2},
\end{eqnarray*}
and that therefore
\begin{eqnarray*}
\p\left( \pd - \E[\pd] \ge t \right) &\le&
\exp\bigg(-\frac{\theta t}{2} + C_7 q^4 e^{\theta/q^2}\bigg) \\ &\le &
e^{-\frac{1}{2}q^2 t \log\left(\frac{t}{C_7 q^2}\right) + q^2 t} =
e^{-\frac{1}{2} q^2 t \log\left(\frac{t}{e^2 C_7 q^2}\right)}.
\end{eqnarray*}
In a similar way, one obtains the bound for the negative tail, and this concludes the proof of
\eqref{eq:lrg4} and the case $d=4$.

Turn now to the case $d=3$. From \eqref{eq:tailbound} and the fact that $V\le 2/q^2$ we get
as above that for some $C_7>0$ we have for any integer $k\ge 4$ that
$$ \E(|V|^k) = \int_0^{2/q^2} ku^{k-1}\p(|V|>u)du \le \frac{C_7}{\vol(B) q^{k-3}}. $$ To get useful
bounds for $k=2$ and $k=3$, observe that
$$ |V| = \left| \frac{1}{|W-x|} - \frac{1}{|W-y|}\right| = \frac{\Big| |W-x|-|W-y| \Big|}{|W-x|\cdot |W-y|}
\le \frac{|x-y|}{|W-x|\cdot |W-y|}. $$
Therefore
\begin{eqnarray*}
\p(|V|>u) &\le& \p\left(|W-x| \le \left(\frac{|x-y|}{u}\right)^{1/2} \right) \\ & &
+ \p\left(|W-y| \le \left(\frac{|x-y|}{u}\right)^{1/2} \right) \le C_8 \frac{|x-y|^{3/2}}{\vol(B) u^{3/2}}.
\end{eqnarray*}
This gives that
\begin{eqnarray*}
\E(|V|^3)& \le& \int_0^{|x-y|^{-1}} 3u^2 \frac{C_8 |x-y|^{3/2}}{\vol(B) u^{3/2}}du \\ & & +
\int_{|x-y|^{-1}}^{2q^{-2}} 3u^2 \frac{2 \kappa_3}{\vol(B) u^3}du \le \frac{C_9 \log(|x-y|/q)}{\vol(B)},
\end{eqnarray*}
and similarly,
\begin{eqnarray*}
 \E(|V|^2) &\le& \int_0^{|x-y|^{-1}} 2u \frac{C_8 |x-y|^{3/2}}{\vol(B) u^{3/2}}du \\ & & +
\int_{|x-y|^{-1}}^{2q^{-2}} 2u \frac{2\kappa_3}{\vol(B) u^3}du \le \frac{C_{10} |x-y|}{\vol(B)}.
\end{eqnarray*}
Combining these bounds and proceeding with the same technique
as above, we deduce that for any $t\ge 0$ and
$\theta \ge0$ we have
\begin{multline*}
\p(\pd - \E[\pd]\ge t) \le \exp\bigg(-\theta t + \frac{C_{10}}{2}\theta^2 |x-y| + \frac{C_9}{6} \theta^3
\log\left(\frac{|x-y|}{q}\right) \\ + C_7 q^3 \sum_{k=4}^\infty \frac{(\theta/q)^k}{k!}\bigg).
\end{multline*}
Take $\theta = \frac{At}{|x-y|}$, where $A$ is a constant such that $0<A<\frac{1}{10(C_7\vee C_9 \vee C_{10}\vee 1)}$.  From the assumptions $2<q<t<|x-y|$, we get 
\begin{eqnarray*}
\theta t &=& \frac{A t^2}{|x-y|}, \\
\frac{C_{10} \theta^2 |x-y|}{2\theta t} &=& \frac{C_{10}A^2 t^2}{2|x-y|\theta t} =\frac{1}{2}C_{10}A < \frac{1}{10},\\
\frac{C_9 \theta^3\log\left(\frac{|x-y|}{q}\right)}6{\theta t} &=& \frac{C_9 A^2 t
\log\left(\frac{|x-y|}{q}\right)}{6|x-y|^2} \le C_9 A^2 < \frac{1}{10},
\end{eqnarray*}
and similarly
$$
\frac{C_7 q^3 \sum_{k=4}^\infty \frac{(\theta/q)^k}{k!}}{\theta t} \le \frac{C_7}{24 \theta t} q^3 \sum_{k=4}^\infty
(\theta/q)^k = \frac{C_7 \theta^3}{24 q t}\frac{1}{1-\theta/q} \le \frac{C_7 A^3}{12} \le \frac{1}{10}
$$
(note that $\theta/q \le 1/2$). Therefore we get
$$ \p(\pd - \E[\pd]\ge t) \le \exp\bigg(-\theta t+\frac{3}{10}\theta t\bigg)
= \exp\left( \frac{-7A t^2}{10 |x-y|}\right). $$
The bound for the negative tail is obtained similarly. This completes the proof of \eqref{eq:lrg44}.
\end{proof}

\begin{cor} \label{largedeviations5}
In dimension $d=4$, there exist constants $C_1, c_2, c_3>0$ such that
for all $x,y\in \R^d$ and $p>q>2$ satisfying $p>3q, |x|,|y|<p/2$, and $|x-y|>3q$, we have that
\begin{multline*}
 \p\left(  \max_{u\in B(x,1), v\in B(y,1)} \Big| \pd\big(u,v\ \big|\ B(0,p) \setminus
(B(u,q)\cup B(v,q)) \big)\Big| > t \right) \\ \le C_1 e^{-c_2 q^2 t \log\left(\frac{c_3 t}{q^2}\right)}
\end{multline*}
for all $t$ that satisfies \eqref{eq:lrgassumption1}.
Similarly, in dimension $d=3$,there exist constants $C_1, c_2, c_3>0$ such that
for all $x,y\in \R^d$, $p>q>2$ and $t>0$ satisfying $p>3q, |x|,|y|<p/2$, and $2<q < t < |x-y|$
 we have that
\begin{multline*}
 \p\left(  \max_{u\in B(x,1), v\in B(y,1)} \Big| \pd\big(u,v\ \big|\ B(0,p) \setminus
(B(u,q)\cup B(v,q)) \big)\Big| > t \right) \\
\le C_1 e^{-\frac{c_2 t^2}{|x-y|}} + C_1 e^{-c_2 q^2
t\log\left(\frac{c_3 t}{q}\right)}.
\end{multline*}
\end{cor}


\begin{proof} This follows from Theorem \ref{largedeviations4} in the same way that
Theorem \ref{largedeviations3} follows from Theorem \ref{largedeviations1}. We omit the proof.
\end{proof}

We also need large deviations estimates for the truncated potential function. This differs from our estimates in dimensions 5 and higher in that the estimates are valid only in a restricted range of the parameters, depending on the dimension.

\begin{thm} \label{largedeviations6}
There exist constants $C_1, c_2, c_3>0$ such that
for all $p>q>0$ we have that
\begin{equation}\label{eq:lrg4second}
 \p\left(  \Big| U\big(0\ \big|\ B(0,p) \setminus B(0,q) \big)\Big| > t \right) \le C_1 e^{-c_2 q^{d-2} t \log\left(\frac{c_3 t}{q^2}\right)}
\end{equation}
for all $t$ above a threshold that depends on $d, p$ and $q$, as follows:
\begin{equation}
\begin{array}{ll} t\ge C_1 q^2\text{ and } \\ t\ge C_1 q^2 \log\left(\frac{C_1 t}{q^2}\right)
  \log\left(\frac{p}{q}\right)\qquad\  & \text{in dimension }d=4; \end{array} \label{eq:lrgassumption3}
\end{equation}
\begin{equation}
\begin{array}{ll}
   t\ge C_1 q^2\text{ and } \\ t\ge C_1 q^2 \left(\log\left(\frac{C_1 t}{q^2}\right)\right)^2
  \log\left(\frac{p}{q}\right)\text{ and } \\
  t\ge C_1 p  q \log\left(\frac{C_1 t}{q^2}\right)
   & \text{in dimension }d=3.\end{array} \label{eq:lrgassumption4}
\end{equation}
\end{thm}

\begin{proof}
Let $d=4$. In the notation of Theorem \ref{largedeviations1}, we now have that \eqref{eq:moment} holds only for $k\ge 3$. For $k=2$ we have
$$ \E|V|^2 = \E|W_{p,q}|^{-4} = \frac{4}{\vol(B)} \int_q^p \kappa_4 t^3 t^{-4} dt = \frac{4\kappa_4}{\vol(B)}
\log\left(\frac{p}{q}\right). $$
Now proceed exactly as in the proof of Theorem \ref{largedeviations4} above.

For $d=3$, we have, still in the notation of Theorem \ref{largedeviations1}, that \eqref{eq:moment} holds only for $k\ge 4$. For $k=2$ we have
$$ \E|V|^2 = \E|W_{p,q}|^{-2} = \frac{3\kappa_3}{\vol(B)}\int_q^p t^2 t^{-2}dt = \frac{3\kappa_3(p-q)}{\vol(B)} \le \frac{3\kappa_3 p}{\vol(B)}, $$
and similarly for $k=3$ we have
$$\E|V|^3 = \frac{3\kappa_3}{\vol(B)} \log\left(\frac{p}{q}\right). $$
Proceeding as in the proofs above, this leads to the inequality
$$ \P\left(U_{p,q}-\E U_{p,q} \ge t \right) \le \exp\left( -\theta t + C_{20} \theta^2 p + C_{30} \theta^3 \log\left( \frac{p}{q}\right) + C_{40}q^3 e^{\theta/q}\right) $$
valid for all $t\ge 0$ and $\theta\ge 0$. Taking $\theta = q\log \left(\frac{t}{C_{40} q^2}\right)$, this easily gives the positive-tail half of \eqref{eq:lrg4second} under the assumptions \eqref{eq:lrgassumption4}. As before the negative-tail half is proved similarly.
\end{proof}

\begin{thm} \label{largedeviations7}
There exist constants $C_1, c_2, c_3, C_4>0$ such that
for all $p>q>0$ we have that, if $t\ge C_4 p^2$, and if the same assumptions
\eqref{eq:lrgassumption3} and \eqref{eq:lrgassumption4} as in Theorem \ref{largedeviations6} hold, then we have
\begin{equation}\label{eq:lrg4other}
 \p\left(  \max_{x\in B\left(0,1\wedge\frac{q}{2}\right)}
 \Big| U\big(x\ \big|\ B(x,p) \setminus B(x,q) \big)\Big| > t \right) \le C_1 e^{-c_2 q^{d-2} t \log\left(\frac{c_3 t}{q^2}\right)}.
\end{equation}
\end{thm}

\begin{proof} This follows from Theorem \ref{largedeviations6} in the same way that
Theorem \ref{largedeviations3} follows from Theorem \ref{largedeviations1}.
\end{proof}

\subsection{Dimension 4 \label{subsectiondim4}}

Let $B, \Delta, r, \rho, s, \eps$ be the same as in Section \ref{subsectiondiscrete}.
We redefine the events $\Omega_1, \Omega_2$, as follows.
\begin{eqnarray*}
\Omega_* &=& \Big\{ \max_{x\in Q(0,2R)}U\big(x\ \big|\ B(x,3R^{1/d})\big) > \frac{B}{2} \Big\}, \\
\Omega_{**} &=& \Big\{ \max_{x \in Q(0,2R)}\Big(\#\text{ of stars in }B(x,3R^{1/d})\Big)
\ge \frac{BR}{10\cdot 100^d} \Big\}, \\
\Omega_1 &=& \Omega_* \cup \Omega_{**}, \\
\Omega_2 &=& \Big\{ \max_{\scriptsize \begin{array}{c}x,y\in Q(0,2R),\\|x-y|\ge R/100\end{array}}
\!\!\!\Big|\pd\Big(x,y\Big|\R^d\!\setminus\! \big(B(x,3R^{1/d})
\cup B(y,3R^{1/d}) \big)\Big)
\Big|\!>\! \frac{B}{2} \Big\}.
\end{eqnarray*}
Let $\Omega_3$ remain the same as in Section \ref{subsectiondiscrete}. The following lemma replaces Lemma \ref{exceptionalomega}.

\begin{lem} \label{exceptionalomega34}  In dimension $4$, for some constants $C,c>0$ we have for all $R>2$ that
\begin{eqnarray}
\p(\Omega_1) &\le& C e^{-c R^{17/9}},  \label{eq:stillholds1} \\
\p(\Omega_2) &\le& C e^{-c R^{11/9}},  \label{eq:exceptest} \\
 \p(\Omega_3) &\le&C e^{-c R^{1+1/100 d^2}} \label{eq:stillholds2}.
\end{eqnarray}
\end{lem}

\begin{proof} First, note that $\Omega_* = \emptyset$ for $R$ sufficiently large, since the only positive contribution to $U\big(x\ \big|\ B(x,3R^{1/d})\big)$ comes from its expected value $(9d\kappa_d/2)
R^{2/d}$ (see \eqref{eq:expectationupq}).

Next, to estimate the probability of $\Omega_{**}$, cover $Q(0,2R+3R^{1/d})$ with $O(R^{d-1})$ balls of radius $6R^{1/d}$ so that for each $x\in Q(0,2R)$, the ball $B(x,3R^{1/d})$ is contained in one of them. In each of these balls, we need to estimate the probability that a Poisson random variable with mean $C\cdot R$ (for some constant $C>0$) is $\ge BR/(10\cdot 100^d)$. By Lemma \ref{poissonlemma}, this probability is $O\left(e^{-c B R \log R}\right)$ for some constant $c>0$.

Finally, the estimate for the probability of $\Omega_2$ follows from Corollary \ref{largedeviations5} in the same way that \eqref{eq:omega2} follows from Theorem \ref{largedeviations3}.
\end{proof}

With these new definitions, the only further change required in the proof of Theorem 3 is the following revised proof of Lemma \ref{manybad}. All the other proofs remain correct as written, with appeals to \eqref{eq:ballfirstonemov} being replaced by using \eqref{eq:lrg4other} instead (one has to verify that the conditions under which \eqref{eq:lrg4other} may be used actually hold, but this is easy).

\begin{proof}[Proof of Lemma \ref{manybad} in dimension 4]
Let $\Gamma:[0,T]\to \R^d$ be a gravitational flow curve that connects $\partial Q(0,R)$ and $\partial Q(0,2R)$, and assume that $\Omega_1\cup \Omega_2$ did not occur. The potential difference $\pd(\Gamma(0),\Gamma(t))$ decreases as a function of $t$. Let
$T_0 = \sup\big\{ t\in [0,T]\,:\,\big|\Gamma(t)-\Gamma(0)\big|\le \frac{R}{100} \big\}$.
Let $W'$ be the set of points $w\in S$ such that $\Gamma_{\big| [T_0,T]}$ intersects both $\partial Q_{\textrm{in}}(w)$
and $\partial Q_{\textrm{out}}(w)$. Since $\Gamma$ connects $\partial Q(0,R)$ and $\partial Q(0,2R)$, we have that $|W'|\ge 99R/100r-2$ (again the $-2$ is to account for boundary effects).

Let
$T_1 = \sup\big\{ t\in [T_0,T]\,:\,\pd(\Gamma(0),\Gamma(t))\ge -2 B \big\}$.
Now, if $x = \Gamma(t)$ for some $t>T_1$, then by the definition of $T_1$ we have
$$
\pd(\Gamma(0),x) \le -2B, $$
and by $\Omega_2^c$ we have
$$
\Big|\pd\Big(\Gamma(0),x\ \Big|\ \R^d\setminus \big(B(\Gamma(0),3R^{1/d})
\cup B(x,3R^{1/d}) \big)\Big)\Big| \le \frac{B}{2}.
$$
Therefore also
$$
\Big|\pd\Big(\Gamma(0),x\ \Big|\ B(\Gamma(0),3R^{1/d})
\cup B(x,3R^{1/d}) \Big)\Big| \le -\frac{3B}{2}.
$$
But
$$
\pd\Big(\Gamma(0),x\ \Big|\ B(\Gamma(0),3R^{1/d})
\cup B(x,3R^{1/d}) \Big) \qquad\qquad\qquad $$

\vspace{-20.0pt}
\begin{eqnarray*}
&=& U\Big(x\ \Big|\ B(x,3R^{1/d})\Big) - U\Big(\Gamma(0)\ \Big|\ B(\Gamma(0),3R^{1/d})\Big)
\\ & &
+ U\Big(x\ \Big|\ B(\Gamma(0),3R^{1/d})\Big) -
U\Big(\Gamma(0)\ \Big|\ B(x,3R^{1/d})\Big),
\end{eqnarray*}
and by $\Omega_1^c$ we have that
$$ -U\Big(\Gamma(0)\ \Big|\ B(\Gamma(0),3R^{1/d})\Big) \ge -\frac{B}{2} $$
and that
$$ \bigg| U\Big(x\ \Big|\ B(\Gamma(0),3R^{1/d})\Big) -
U\Big(\Gamma(0)\ \Big|\ B(x,3R^{1/d})\Big)\bigg| \qquad\qquad\qquad\qquad $$
$$ \qquad\qquad \le \frac{BR}{10\cdot 100^d} \cdot \frac{1}{(R/100)^{d-2}} \le \frac{B}{10}. $$
Therefore we get that, for $x=\Gamma(t), t>T_1$, we have
$$ U\Big(x\ \Big|\ B(x,3R^{1/d})\Big) \le -\frac{9B}{10} < -\frac{B}{2}. $$
By the above, it follows  that all the points $w\in S$ for which $\Gamma_{\big| [T_1,T]}$ intersects both $\partial Q_{\textrm{in}}(w)$ and
$\partial Q_{\textrm{out}}(w)$ are bad, since they satisfy condition (1) in the definition.
If there are $R/10^d r$ such points, we are done, since the set of such $w$ is connected and \emph{a fortiori} connectible. If this is not so, denote by $W''$ the set of those $w \in W'$ for which
$\Gamma_{\big| [T_0,T_1]}$ intersects both $\partial Q_{\textrm{in}}(w)$ and $\partial Q_{\textrm{out}}(w)$. The family $W''$ is
a connected set, and we have $|W''| \ge |W'| - R/4r - 1 \ge 74R/100r-4 > R/2r$.
As in Section \ref{subsectiondiscrete}, replace $W''$ by a subset $W'''\subset W''$ such that $|W'''|\ge |W''|/5^d$ and all the interiors
of $(Q_{\text{out}}(w))_{w\in W'''}$ are disjoint. Repeating the same argument as in Section \ref{subsectiondiscrete}, we get that the set $W$ of bad $w\in W'''$ is connectible
and contains $\ge R/10^d r$ points.
%
\end{proof}

\subsection{Dimension 3 \label{subsectiondim3}}

Let $d=3$. All the foregoing discussion for dimension 4 remains valid, except the estimate
\eqref{eq:exceptest}. In dimension 3 we only get the weaker estimate
$$ \p(\Omega_2) \le Ce^{-cR^{7/9}} $$
for some constants $C,c>0$. Thus, while all the elements of the proof still function, what we
actually proved was an upper bound for $\p(E_R)$ which is of the form $Ce^{-cR^{7/9}} $.

To get the better bound stated in Theorem \ref{thm3}, we modify the value of the parameters.
Here are the new values:
\begin{eqnarray*}
B&=& \frac{R}{(\log R)^\beta}, \\
r&=& (\log R)^{1/3}\log\log R, \\
\rho&=& \frac{1}{(\log R)^\gamma}, \\
s&=& \frac{1}{(\log R)^\delta}, \\
\varepsilon&=& \frac{\rho}{s^3}\log \log R, \\
\lambda &=& \sqrt{\log\log R}.
\end{eqnarray*}
Here $\beta, \gamma, \delta$ are positive constants. The proof in Section \ref{sectionproofthm3}, together with the adjustments of of Subsection \ref{subsectiondim4},
will work almost verbatim with these modified parameters, provided several conditions are met:

\begin{itemize}
\item $ \lambda < c \frac{r}{(\log R)^{1/3}}$ for some constant $c>0$ and all sufficiently large $R$. This is required when using
Theorem \ref{joint_force_density_thm}, and holds with our choice of parameters.
\item $\varepsilon>C \frac{B}{R}$ for some constant $C>0$. This is used in the proof of Lemma
\ref{pwbound}, when we deduce from $|F(y)|<CB/R$ that in fact $|F(y)|<\eps/8$.
It will hold if $\gamma<\beta+3\delta$, and in particular if $\gamma <\beta$.
\item
$\left(\frac{r}{\rho}\right)^{d-1} \eps^d \lambda^{d^2-d} << 1$.
This is
required when using Markov's inequality to ensure that the probability per site $w\in W$
to have $\ge \frac{r}{2\cdot 7^d\rho}$ points $x\in T_w$ which are black and not $3$-crowded
(Lemma \ref{pwbound}) is $<<1$. This condition holds for any
$\gamma>\frac{2}{3}+9\delta$, and in particular for any
$\gamma>\frac{2}{3}$ if $\delta$ is sufficiently small (as a function
of $\gamma$).
\item $\frac{s^d \eps}{\rho} >> 1$. This is required in the proof of
Lemma \ref{pwbound} to make sure that $\p(A_j') \le \frac{1}{R^d}$.
\end{itemize}
With this choice of parameters, following the steps of the proof in Section \ref{sectionproofthm3}
together with the changes outlined in Subsection \ref{subsectiondim4}, we get that
in Lemma \ref{exceptionalomega34} the estimates \eqref{eq:stillholds1} and \eqref{eq:stillholds2}
still hold. The estimate \eqref{eq:exceptest} is replaced by the following estimate, whose proof again uses Corollary \ref{largedeviations5}.
\begin{equation} \label{eq:replaced}
\p(\Omega_2) \le Ce^{-c \frac{R}{(\log R)^{2\beta}}}.
\end{equation}
Lemma \ref{pwbound} will be weakened to the following lemma, whose
proof is a repetition of the same steps with the new parameter values.

\begin{lem} \label{pwboundnew} For some constants $C,c>0$, we have
$$ \p(\Omega_{7,W}) \le C e^{-c |W| \log\log R}. $$
\end{lem}

As a result, Theorem \ref{probbad} will be weakened to the following theorem.

\begin{thm} \label{probbadnew} There exist constants $C,c>0$ such that
for any family $W \subseteq S$ we have
\begin{equation} \label{eq:probfamilybad2}
\p\bigg(
\Omega_3^c \cap \Big\{ \text{all }w\in W\text{ are bad}\Big\}
\bigg) \le C e^{-c |W| \log \log R}.
\end{equation}
\end{thm}

Therefore, the bound that we get for $\p(E_R)$ will be weakened to
$$ \p(E_R) \le \p(\Omega_1) +\p(\Omega_2) +\p(\Omega_3) + e^{-c \frac{B}{r} \log\log R}
\le Ce^{-c\frac{R}{(\log R)^{2\beta}}} + Ce^{-c \frac{R}{(\log R)^{1/3}}}. $$
With the constraints $\beta + 3\delta > \gamma > 2/3 + 9\delta$,
the best that one can do is to take $\beta$ slightly bigger than $2/3$. This gives that for all $\alpha>4/3$ we have
$$ \p(E_R) \le Ce^{-c \frac{R}{(\log R)^\alpha}} $$
for some constants $c,C>0$ depending on $\alpha$ and all $R>2$. This was the claim of Theorem
\ref{thm3} in dimension 3.
\qed

\bigskip \noindent
Sourav Chatterjee \\
Department of Statistics \\
367 Evans Hall \\
The University of California \\
Berkeley, CA 94720-3860, USA \\
\texttt{sourav@stat.berkeley.edu}

\bigskip \noindent
Ron Peled\\
Department of Statistics \\
367 Evans Hall \\
The University of California \\
Berkeley, CA 94720-3860, USA \\
\texttt{peledron@stat.berkeley.edu}

\bigskip \noindent
Yuval Peres \\
Microsoft Research \\
One Microsoft way \\
Redmond, WA 98052-6399, USA \\
\texttt{peres@stat.berkeley.edu}

\bigskip\noindent
Dan Romik \\
Bell Laboratories \\
Fundamental and Industrial Mathematics Research Department \\
Room 2C-379 \\
600 Mountain Ave \\
Murray Hill, NJ 07974, USA \\
\texttt{romik@bell-labs.com}

\end{document}